\documentclass{amsart}
\usepackage[utf8]{inputenc}
\usepackage{amsmath}
\usepackage{amssymb}
\usepackage{amsfonts}
\usepackage{amscd}
\usepackage{amsthm}
\usepackage{xypic}
\usepackage{enumerate}
\usepackage{fancyhdr}

\newcommand{\Z}{\mathbb{Z}}
\newcommand{\Q}{\mathbb{Q}}

\newcommand{\F}{\mathbb{F}}
\newcommand{\K}{\mathbb{K}}

\newcommand{\mmP}{\mathcal{P}}

\newcommand{\mmO}{\mathcal{O}}
\newcommand{\mmC}{\mathcal{C}}

\newcommand{\mmF}{\mathcal{F}}

\newcommand{\mmK}{\mathcal{K}}

\newcommand{\mmD}{\mathcal{D}}

\newcommand{\uh}{\underline{\Hom}}

\DeclareMathOperator*{\holim}{holim} \DeclareMathOperator*{\pr}{pr}
 
\DeclareMathOperator*{\Hom}{Hom} \DeclareMathOperator*{\Ho}{Ho}

\DeclareMathOperator*{\Map}{Map}
\DeclareMathOperator*{\Sub}{Sub}\DeclareMathOperator*{\Obj}{Obj}
\DeclareMathOperator*{\Vect}{Vect}
\DeclareMathOperator*{\sk}{sk}\DeclareMathOperator*{\Zar}{Zar}
\DeclareMathOperator*{\ZAR}{ZAR}\DeclareMathOperator*{\ch}{ch}

\numberwithin{equation}{section}
\newtheorem{theo}[equation]{Theorem}
\newtheorem{prop}[equation]{Proposition}
\newtheorem{cor}[equation]{Corollary}
\newtheorem{lema}[equation]{Lemma}

\theoremstyle{definition}
\newtheorem{df}[equation]{Definition}
\newtheorem{obs}[equation]{Remark}
\newtheorem{ex}[equation]{Example}

\newenvironment{enumerate*}[1][{}]{\begin{itemize}\renewcommand{\labelitemi}{#1}}{\end{itemize}}

\title{On uniqueness of characteristic classes}
\author{Elisenda Feliu}
\email{efeliu@ub.edu}
\address{Gran Via de les Corts Catalanes, 585, 08007 Barcelona (Spain)}

\newif\ifprivate
\privatetrue

\begin{document}

\begin{abstract}
We give an axiomatic characterization of maps from algebraic K-theory. The results
apply to a class of maps from algebraic K-theory to any suitable cohomology theory or to algebraic K-theory, which includes all group morphisms. In particular, we obtain comparison theorems for the Chern
character and Chern classes and for the Lambda and Adams operations on higher algebraic K-theory. We show that the Adams operations defined by Grayson agree with the ones defined by Gillet and Soul\'e.
\end{abstract}

\maketitle

\section*{Introduction}
In this paper we address the problem of comparing maps from the algebraic $K$-groups of a scheme to either algebraic $K$-groups or
suitable cohomology theories. This type of questions often arise when one constructs a map that is supposed to induce, in some particular cohomology theory, a specific regulator or known map, and one needs to show that the map is indeed the expected one. Examples of these situations are found in the construction of the Beilinson regulator given by Burgos and Wang in \cite{Burgos1} or in the definition of the Adams operations given by Grayson in \cite{Grayson1}.

The different nature of each construction usually makes direct comparisons not an available option and one is forced to turn to  theoretical tricks. In this work, we identify sufficient conditions for two maps to agree, thus obtaining an axiomatic characterization of maps from $K$-theory.

\par
In abstract, the results apply to a class of maps, named \emph{weakly additive}. All group morphisms induced by a map of sheaves are in this class, but are not the only ones.  As a main consequence, we give a characterization of the
\emph{Adams}  and \emph{lambda operations} on higher $K$-theory and of the \emph{Chern character} and \emph{Chern classes} on a suitable cohomology theory (see section \ref{uniqueness4}).

In particular, we show that the Adams operations defined by Grayson in \cite{Grayson1} agree with the ones defined by Gillet and Soul{\'e} in \cite{GilletSouleFiltrations}, for all noetherian schemes of finite Krull dimension. This implies that for this class of schemes, the operations defined by Grayson satisfy the usual identities of a lambda ring. This is an original result.

The second specific application of this work is a proof that the regulator defined by Burgos and Wang in \cite{Burgos1} is the Beilinson regulator. The proof provided here is simpler than the one given in loc. cit., where delooping in $K$-theory was required.

The results of this paper are further exploited in the paper under preparation by the author (\cite{Feliu-Adams},\cite{FeliuThesis}), where an explicit chain morphism representing the Adams operations on higher algebraic $K$-theory with rational coefficients is constructed.
Furthermore, in \cite{FeliuThesis}, \cite{FeliuChow}, Burgos and the author defined a morphism in the derived category of complexes from a chain complex computing higher algebraic Chow groups to Deligne-Beilinson cohomology. A slight modification of the tools developed here allows one to prove that this morphism induces the Beilinson regulator. This result ultimately implies that the regulator defined by Goncharov in \cite{Goncharov} induces the Beilinson regulator as well (work under preparation by Burgos, Feliu and Takeda).

\par The techniques used in this paper rely on the generalized cohomology theory described by Gillet and Soul{\'e} in
\cite{GilletSouleFiltrations}. Roughly speaking, the idea is that any good enough map from $K$-theory to $K$-theory or to a
cohomology theory is characterized by its behavior over the $K$-groups of the simplicial classifying scheme $B_{\cdot}GL_N$.

More explicitly, let $\mathbf{C}$ be the big Zariski site over a noetherian finite dimensional scheme $S$. For any integers $N,k\geq 0$ denote by $B_{\cdot}GL_{N/S}$ the simplicial scheme $B_{\cdot}GL_N\times_{\Z} S$ and let $Gr_S(N,k)$ be the Grassmanian
scheme over $S$. If $S_{\cdot}\mmP$ denotes the Waldhausen simplicial sheaf computing algebraic $K$-theory, for every scheme $X$  in $\mathbf{C}$ and for all $m\geq 0$, we have $K_m(X) \cong \pi_{m+1}(S_{\cdot}\mmP(X))$. The definition of algebraic $K$-groups can be extended to simplicial schemes over $\mathbf{C}$ and every sheaf map (in the homotopy category of simplicial sheaves) $S_{\cdot}\mmP \rightarrow S_{\cdot}\mmP$ induces a map $K_m(Y_{\cdot}) \rightarrow K_m(Y_{\cdot})$ for every simplicial scheme $Y_{\cdot}$.

Let $\F_{\cdot}$ be a
simplicial sheaf over $\mathbf{C}$. The two main consequences of our uniqueness theorem are the following.
\emph{\begin{enumerate}[(i)] \item Assume that $\F_{\cdot}$ is
weakly equivalent to $S_{\cdot}\mmP$, and let
$\Phi,\Phi':S_{\cdot}\mmP\rightarrow \F_{\cdot}$ be two H-space
maps. Then, the morphisms
$$\Phi,\Phi': K_m(X) \rightarrow K_m(X),\qquad m\geq 0$$ agree for all schemes $X$ in $\mathbf{C}$,
if they agree over $K_0(B_{\cdot}GL_{N/S})$, for all $N\geq 1$.
\item Assume that $\F_{\cdot}$ is weakly equivalent to
$\mmK_{\cdot}(\mmF(*))$, where $\mmK_{\cdot}(\cdot)$ is the sheaf
version of the Dold-Puppe functor, and $\mmF(*)$ is a graded sheaf
giving a twisted duality cohomology theory in the sense of Gillet (\cite{Gillet}). Let $\Phi,\Phi':S_{\cdot}\mmP\rightarrow
\F_{\cdot}$ be two H-space maps. Then, the morphisms
$$\Phi,\Phi': K_m(X) \rightarrow H^{*}(X,\mmF(*)),\qquad m\geq 0$$ agree for all schemes $X$ in $\mathbf{C}$,
if they agree over $K_0(Gr_S(N,k))$, for all $N,k$.
\end{enumerate}}

It follows from \emph{(i)} that there is a unique way to extend the Adams operations from the Grothendieck of vector bundles over a scheme $X$ to higher $K$-theory by means of a sheaf map (in the homotopy category of simplicial sheaves). Analogously, the result \emph{(ii)} implies that there is a unique way to extend the Chern character of vector bundles over a scheme $X$ to higher degrees by means of a sheaf map (in the homotopy category of simplicial sheaves).

The paper is organized as follows. The first two sections are
dedicated to review part of the theory developed by Gillet and
Soul{\'e} in \cite{GilletSouleFiltrations}. More concretely, in
Section 1 we recall the main concepts about the homotopy theory of
simplicial sheaves and generalized cohomology theories. In Section
2 we explain how $K$-theory can be given in this setting. We
introduce the simplicial sheaves $\K^N_{\cdot}=\Z\times
\Z_{\infty}B_{\cdot}GL_N$ and $\K_{\cdot}=\Z\times
\Z_{\infty}B_{\cdot}GL$.

In Section 3, we consider compatible systems of maps
$\Phi_{N}:\K^N_{\cdot}\rightarrow \F_{\cdot}$, for $N\geq 1$ and
$\F_{\cdot}$ a simplicial sheaf. We introduce the class of
\emph{weakly additive} system of maps, which are the ones
characterized in this paper. Roughly speaking, they are the
systems for which all the information can be obtained separately
from the composition $\Z_{\infty}B_{\cdot}GL_N\hookrightarrow
\K^N_{\cdot}\rightarrow \F_{\cdot}$ and from the composition
$\Z\hookrightarrow \K^N_{\cdot}\rightarrow \F_{\cdot}$. They are
named \emph{weakly additive} due to the fact that when
$\F_{\cdot}$ is an $H$-space, they are the maps given by the sum
of this two mentioned compositions. The main example are the
systems of the type $\K^{N}_{\cdot}\hookrightarrow \K_{\cdot}
\xrightarrow{\Phi} \F_{\cdot}$, inducing group morphisms on
cohomology. We discuss then the comparison of two different weakly
additive systems of maps. We end this section by applying the
general discussion  to generalized cohomology theories on a
Zariski site.

In the last two sections we develop the application of the characterization
results to $K$-theory and to cohomology theories. Section 4 is devoted to maps from $K$-theory to $K$-theory, concretely to the
Adams and lambda operations on higher algebraic $K$-theory. A characterization of these operations is given and the comparison to Grayson Adams operations is provided.
In section 5, we consider maps from  the $K$-groups to suitable sheaf cohomologies. We
give a characterization of the Chern character and of the Chern classes for the
higher algebraic $K$-groups of a scheme.

\emph{Acknowledgement.} I would like to thank Jos{\'e} Ignacio Burgos Gil for
his dedication in the elaboration of this work. I am grateful to Christophe
Soul{\'e} for  a  useful discussion about homotopy theory of sheaves.

\section{The homotopy category of simplicial sheaves}

We review here the main definitions and properties of homotopy
theory of simplicial sheaves. For more details about this topic
see \cite{GilletSouleFiltrations}. For  general facts and
definitions about model categories we refer for instance to
\cite{Hirschhorn}.

Let $\mathbf{C}$ be a site and let $\mathbf{T}=T(\mathbf{C})$ be the
(Grothendieck) topos of sheaves on $\mathbf{C}$. We will always suppose that
$\mathbf{T}$ has enough points (see \cite{SGA4}, $\S$ IV 6.4.1).

Let $\mathbf{sT}$ be the category of simplicial objects in
$\mathbf{T}$. One identifies $\mathbf{sT}$ with the category of
sheaves of simplicial sets on $\mathbf{C}$.  An object of
$\mathbf{sT}$ is called a \emph{space}.

\subsection{Structure of simplicial model category}\label{def} The category $\mathbf{sT}$
is endowed with a  structure of simplicial model category in the sense of Quillen
 \cite{Quillen0}. This result is due to Joyal;
a proof of it can be found in \cite{Jarsimppre}, corollary. 2.7. Here we recall
the definitions that give a simplicial model structure to $\mathbf{sT}$.

The structure of \emph{model category} of $\mathbf{sT}$ is given
as follows. Let $X_{\cdot}$ be a space in $\mathbf{sT}$. One
defines $\pi_0(X_{\cdot})$ to be the sheaf associated to the
presheaf
$$U\mapsto \pi_0(X_{\cdot}(U)),\quad \textrm{for } U\in \Obj(\mathbf{C}).$$

Let $\mathbf{C}|U$ be the site of objects over $U$ as described in
\cite{SGA4}, $\S$ III 5.1, and let $\mathbf{T}|U$ denote the
corresponding topos. For every object $X_{\cdot}$ in
$\mathbf{sT}$, let $X_{\cdot}|U$ be the restriction of $X_{\cdot}$
to $\mathbf{sT}|U$. Then, for every $U\in \Obj(\mathbf{C})$, $x\in
X_0(U)$ a vertex of the simplicial set $X_{\cdot}(U)$, and every
integer $n>0$, one defines $\pi_n(X_{\cdot}|U,x)$ to be the sheaf
associated to the presheaf
$$V \mapsto \pi_n(X_{\cdot}(V),x),\quad \textrm{for } V\in \Obj(\mathbf{C}|U).$$

Let $X_{\cdot},Y_{\cdot}$ be two spaces and let
$f:X_{\cdot}\rightarrow Y_{\cdot}$ be a map.
\begin{enumerate}[(i)]
\item The map $f$ is called a \emph{weak equivalence} if the
induced map $f_*:\pi_0(X_{\cdot})\rightarrow \pi_0(Y_{\cdot})$ is
an isomorphism and, for all $n>0$, $U\in \Obj(\mathbf{C})$ and
$x\in X_0(U)$, the natural maps
$$f_*: \pi_n(X_{\cdot}|U,x)\rightarrow \pi_n(Y_{\cdot}|U,f(x))$$ are isomorphisms.
\item The map $f$ is called a \emph{cofibration}
 if for every $U\in
\Obj(\mathbf{C})$, the induced map
$$f(U): X_{\cdot}(U)\rightarrow Y_{\cdot}(U)$$
 is a cofibration of simplicial
sets, i.e. it is a monomorphism. \item The map $f$ is called a
\emph{fibration}   if it has the right
lifting property with respect to trivial cofibrations.
\end{enumerate}
Observe that since the only map $\emptyset\rightarrow X_{\cdot}$
is always a monomorphism, all objects $X_{\cdot}$ in $\mathbf{sT}$
are cofibrant.

Let $\mathbf{SSets}$ denote the category of simplicial sets. The
structure of \emph{simplicial category} of $\mathbf{sT}$ is given
by the following definitions:
\begin{enumerate}[(i)]
\item There is a functor $\mathbf{SSets} \rightarrow \mathbf{sT}$,
which sends every simplicial set $K_{\cdot}$ to the sheafification
of the constant presheaf that takes the value $K_{\cdot}$ for
every $U$ in $\mathbf{C}$. \item For every space $X_{\cdot}$ and
every simplicial set $K_{\cdot}$, the direct product
$X_{\cdot}\times K_{\cdot}$ in $\mathbf{sT}$ is the simplicial set
given by
$$[n]\mapsto \coprod_{\sigma\in K_n}X_n,$$ and induced face and
degeneracy maps.
 \item Let $X_{\cdot}$, $Y_{\cdot}$ be
two spaces and let $\Delta^n_{\cdot}$ be the standard $n$-simplex
in $\mathbf{SSets}$. The simplicial set $\uh(X_{\cdot},Y_{\cdot})$
is the functor
$$[n]\mapsto \Hom\nolimits_{\mathbf{sT}}(X_{\cdot}\times \Delta^n_{\cdot},Y_{\cdot}).$$
\end{enumerate}

Note that by definition, a map is a cofibration of spaces if and only if it is a
section-wise cofibration of simplicial sets. For fibrations and weak
equivalences, this is not always true. However, it follows from the definition
that   a section-wise weak equivalence is a weak equivalence of spaces.

\subsection{Fibrant resolutions and the homotopy category}\label{fact} Let $\mmC$ be any model category and
$f:X\rightarrow Y$
a map. By definition,  there exist two  factorizations $(\alpha,\beta)$ of $f$,
\begin{enumerate}
\item $f=\beta\circ \alpha$, with $\alpha$ a cofibration and $\beta$ is a
trivial fibration, \item $f=\delta\circ \gamma$, where $\gamma$ is a trivial
cofibration and $\delta$ is a fibration.
\end{enumerate}

A \emph{fibrant resolution} of an object $X$, is a fibrant object $X^{\sim}$
together with a trivial cofibration $X\xrightarrow{\sigma} X^{\sim}$.
\emph{Cofibrant resolutions} are defined dually.  Its existence is guaranteed
by the existence of the factorizations above.

One can form the \emph{homotopy category} $\Ho (\mathbf{sT})$, associated to
$\mathbf{sT}$, by formally inverting the weak equivalences. For any two spaces
$X,Y$, one denotes by $[X,Y]$ the set of maps between $X$ and $Y$ in this
category.

If $Y\rightarrow Y^{\sim}$ is any fibrant resolution of $Y$, then
$$[X,Y]=\pi_0 \Hom(X,Y^{\sim}).$$

More generally, suppose that $Y\rightarrow \hat{Y}$ is a weak equivalence (not necessarily also a cofibration) and
$\hat{Y}$ is fibrant. Then, if $Y^{\sim}$ is any fibrant resolution, there
exists a weak equivalence $Y^{\sim}\rightarrow \hat{Y}$. Therefore, by
\cite{Hirschhorn}, 9.5.12,
$$[X,Y]=\pi_0 \Hom(X,\hat{Y}).$$

Consider $X,Y\in \Obj(\mmC)$ and $f:X\rightarrow Y$ a morphism. Suppose that
$Y$ is fibrant and let $X^{\sim}$ be a fibrant resolution of $X$. Then, $f$
factors uniquely (up to homotopy) through $X^{\sim}$, i.e. there exists a map
in $\mmC$, $f^{\sim}:X^{\sim}\rightarrow Y$, unique up to homotopy under $X$,
such that the following diagram is commutative
$$
\xymatrix{ X \ar[r]^f \ar[d]_{\sigma}^{\sim} & Y \\   X^{\sim}
\ar[ur]_{f^{\sim}}. }
$$
See \cite{Hirschhorn}, 8.1.6 for a proof. Therefore, there is a map
$$[X,Y] \rightarrow [X^{\sim},Y]$$
obtained by the factorizations.

\subsection{The category of simplicial presheaves} Let $\mathbf{sPre(C)}$ be the
category of simplicial presheaves on $\mathbf{C}$, i.e. the
  category of functors $\mathbf{C}^{op}\rightarrow
\mathbf{SSets}$. Then, one defines:
\begin{enumerate*}[$\rhd$] \item \emph{weak equivalences} and \emph{cofibrations} of
simplicial presheaves exactly as for simplicial sheaves, and,
\item \emph{fibrations} to be the maps satisfying the right
lifting property with respect to trivial cofibrations.
\end{enumerate*}
As shown by Jardine in \cite{Jarsimppre}, these definitions equip
$\mathbf{sPre(C)}$ with a model category structure. The sheafification functor
$$\mathbf{sPre(C)} \xrightarrow{ ({\cdot})^s} \mathbf{sT} $$
induces an equivalence between the respective homotopy categories,
sending weak equivalences to weak equivalences. Moreover, the
natural map $X_{\cdot}\rightarrow X_{\cdot}^s$ is a weak
equivalence of simplicial presheaves (\cite{Jarsimppre}, Lemma
2.6).

\subsection{The category of pointed simplicial sheaves} Let  $\mathbf{sT_*}$ denote the category of pointed simplicial sheaves with morphisms being the maps preserving the base points. The definitions and results stated above can be translated into this category by considering the analogous definitions for pointed objects
(see \cite{HuberWild} \S B.1 for details).

Let $*$ be the base point of $\mathbf{sT_*}$ (i.e., the final and initial object). Then, if $X_{\cdot}$
is a simplicial sheaf,  one considers its associated pointed
object to be $X_{{\cdot}+}=X_{\cdot}\sqcup *$.

\subsection{Generalized cohomology theories}
Let $X_{\cdot}$ be any space in $\mathbf{sT_*}$. The
\emph{suspension} of $X_{\cdot}$, $S\wedge X_{\cdot}$, is defined
to be the space $X_{\cdot}\wedge \Delta^1/\sim$, where $\sim$ is
the equivalence relation generated by $(x,0)\sim (x,1)$ and where
$\wedge$ is the pointed product. The \emph{loop space functor}
 $\Omega$ is the
right adjoint functor of $S$ in the homotopy category.

Let $A_{\cdot}$ be any  space in $\mathbf{sT}_*$. For every  space
$X_{\cdot}$ in $\mathbf{sT}_*$, one defines the \emph{cohomology
of $X_{\cdot}$ with coefficients in $A_{\cdot}$} as
$$H^{-m}(X_{\cdot},A_{\cdot})=[S^m\wedge X_{\cdot},A_{\cdot}],\quad m\geq 0.$$
This is a pointed set for $m=0$, a group for $m>0$ and an abelian group for
$m>1$.

An infinite loop spectrum $A^*_{\cdot}$ is a collection of spaces
$\{A^i_{\cdot}\}_{i\geq 0}$, together with given  weak equivalences
$A^i_{\cdot}\xrightarrow{\sim} \Omega A^{i+1}_{\cdot}$. \emph{The
cohomology with coefficients in the spectrum $A^*_{\cdot}$} is
defined as
$$H^{n-m}(X_{\cdot},A^*_{\cdot})=[S^m\wedge X_{\cdot},A^n_{\cdot}],\quad m,n\geq 0.$$
Due to the adjointness  relation between the loop space functor and the suspension, these sets depend
only on the difference $n-m$. Therefore, all of them are abelian groups.

Let $A_{\cdot}$ be a simplicial sheaf and assume that there is an
infinite loop spectrum $A^*_{\cdot}$ with $A^0_{\cdot}=A_{\cdot}$.
Then the cohomology groups with coefficients in $A_{\cdot}$ are
also defined with positive indices, with respect to this infinite
loop spectrum. By abuse of notation, when there is no source of
confusion, we will write $H^{m}(X,A_{\cdot})$, for the generalized
cohomology with positive indices, instead of writing
$H^m(X,A_{\cdot}^*)$.

When $X_{\cdot}$ is a non-pointed space in $\mathbf{sT}$, we
define
$$H^{-m}(X_{\cdot},A_{\cdot})=[S^m\wedge X_{{\cdot}+},A_{\cdot}],\quad m\geq 0.$$

\subsection{Induced morphisms} Let $A_{\cdot},B_{\cdot}$ be two pointed spaces. Every element $f\in
[A_{\cdot},B_{\cdot}]$ induces functorial maps
$$ [X_{\cdot},A_{\cdot}]  \xrightarrow{f_*}  [X_{\cdot},B_{\cdot}] \quad \textrm{and}\quad
 [B_{\cdot},X_{\cdot}]  \xrightarrow{f^*}
[A_{\cdot},X_{\cdot}],$$ for every space $X_{\cdot}$. Therefore,
there are induced maps between the generalized cohomology groups
$$H^{-*}(X_{\cdot},A_{\cdot})\xrightarrow{f_*} H^{-*}(X_{\cdot},B_{\cdot})
\quad \textrm{and}\quad
H^{-*}(B_{\cdot},X_{\cdot})\xrightarrow{f^*}
H^{-*}(A_{\cdot},X_{\cdot}).$$

Using simplicial resolutions, these maps can be described as
follows.  If $B^{\sim}_{\cdot}$ is any fibrant resolution of
$B_{\cdot}$, then $f$ is given by a homotopy class of maps
$A_{\cdot}\rightarrow B_{\cdot}^{\sim}$. This map factorizes,
uniquely up to homotopy, through a fibrant resolution of
$A_{\cdot}$, $A_{\cdot}^{\sim}$. Therefore there is a map
$$f^{\sim}:A_{\cdot}^{\sim}\rightarrow B^{\sim}_{\cdot}$$
which induces, for every $m\geq 0$, a map
$$H^{-m}(X_{\cdot},A_{\cdot})=
\pi_0\Hom(S^m\wedge X_{\cdot},A_{\cdot}^{\sim})\rightarrow
\pi_0\Hom(S^m\wedge
X_{\cdot},B^{\sim}_{\cdot})=H^{-m}(X_{\cdot},B_{\cdot}).$$ The
description of $f^*$ is analogous.

\subsection{Zariski topos }\label{zar}
By the \emph{big Zariski site},  $\ZAR$,
we refer to the category of all noetherian schemes of finite Krull
dimension, equipped with the Zariski topology.

Given any scheme $X$, one can consider the category formed by the
inclusion maps $V\rightarrow U$ with $U$ and $V$ open subsets of
$X$ and then define the covers of $U\subseteq X$ to be the open
covers of $U$. This is called the \emph{small Zariski site} of $X$,
$\Zar(X)$. By the \emph{big Zariski site }of $X$, $\ZAR(X)$, we
mean the category of all schemes of finite type over $X$ equipped
with the Zariski topology.

The corresponding topos are named the small or big Zariski topos (over $X$)
respectively.

Generally, one also considers subsites of the big and small
Zariski sites. For instance, the site of all noetherian schemes of
finite Krull dimension which are also smooth (regular,
quasi-projective or projective resp.) is a subsite of $\ZAR$. Similar
subsites can be defined in $\ZAR(X)$ and $\Zar(X)$, depending on
the properties of $X$.

At any of these sites, one associates to every scheme $X$ in the
underlying category $\mathbf{C}$,  the constant pointed simplicial
sheaf
$$U\mapsto \Map\nolimits_{\mathbf{C}}(U,X)\cup \{*\}, \qquad  U\in \Obj \mathbf{C}.$$
 This simplicial sheaf is also denoted by $X$. For any
simplicial sheaf $\F_{\cdot}$ and any scheme $X$ in $\mathbf{C}$,
the equality of simplicial sets
$$\F_{\cdot}(X)= \Hom(X_{\cdot},\F_{\cdot})$$
is satisfied.

\begin{df}
A space $X_{\cdot}$ is said to be \emph{constructed from schemes}
 if, for all $n\geq 0$,
$X_n$ is representable by a scheme in the site plus a disjoint
base point. If $P$ is a property of schemes, one says that
$X_{\cdot}$ satisfies the property $P$, if this is the case for the scheme parts of the components.

\end{df}
Any simplicial scheme gives rise to a space constructed from schemes, but the
converse is not true (see \cite{HuberWild}, {\S}B.1).

If $X_{\cdot}$ is a space constructed from schemes, we can write
$X_n=* \coprod X'_n$, with $X'_n$ a scheme. For every pointed
simplicial sheaf $\F_{\cdot}$ in $\mathbf{sT}_*$, set
$\F_{\cdot}(X_n):=\F_{\cdot}(X'_n).$ Then, one defines
$$\F_{\cdot}(X_{\cdot})= \holim_{\substack{\longleftarrow \\ n}} \F_{\cdot}(X_n), $$
where $\holim$ is the homotopy limit functor defined in
\cite{bousfieldkan}.

\subsection{Pseudo-flasque presheaves}
We fix $\mathbf{T}$ to be the topos associated to any Zariski site
$\mathbf{C}$ as in the previous section. The next definition is at
the end of $\S 2$ in \cite{BrownGestern}.
\begin{df}\label{BG1}
Let $\F_{\cdot}$ be a pointed simplicial presheaf on $\mathbf{C}$.
It is called a \emph{pseudo-flasque} presheaf, if the following
two conditions hold:
\begin{enumerate}[(i)]
\item $\F_{\cdot}(\emptyset)=0$. \item For every pair of open
subsets $U,V$ of some scheme $X$, the square
$$\xymatrix{ \F_{\cdot}(U\cap V) \ar[r] \ar[d] & \F_{\cdot}(V) \ar[d] \\ \F_{\cdot}(U) \ar[r] & \F_{\cdot}(U\cup V)} $$
is homotopy  cartesian.
\end{enumerate}
\end{df}
 A pseudo-flasque presheaf $\F_{\cdot}$ satisfies the Mayer-Vietoris property, i.e.
for any scheme $X$ in the site $\mathbf{C}$ and any two open
subsets $U,V$ of $X$, there is a long exact sequence
{\small $$\cdots \rightarrow H^i(\F_{\cdot}(U\cup V)) \rightarrow H^i(\F_{\cdot}(U)\oplus
\F_{\cdot}(V)) \rightarrow H^i(\F_{\cdot}(U\cap V)) \rightarrow
 H^{i+1}(\F_{\cdot}(U\cup V)) \rightarrow \cdots $$ }

 The importance of pseudo-flasque presheaves relies on the following
 proposition, due to Brown and Gestern (see \cite{BrownGestern}, Theorem 4).
 \begin{prop}\label{BG}
Let $\F_{\cdot}$ be a pseudo-flasque presheaf. For every scheme
$X$ in $\mathbf{C}$, the natural map
$$\pi_i(\F_{\cdot}(X)) \rightarrow H^{-i}(X,\F_{\cdot}^s) $$
is an isomorphism.
 \end{prop}
 Observe that this proposition is already true for any fibrant space. In fact,
 any fibrant space is pseudo-flasque.

\section{K-theory as a generalized cohomology}\label{K}

 Let $X_{\cdot}$ be a space such that its $0$-skeleton is reduced to one point.
 One defines $\Z_{\infty}X_{\cdot}$ to be the sheaf associated to
the presheaf
$$U\mapsto \Z_{\infty}X_{\cdot}(U),$$
the functor $\Z_{\infty}$ being the Bousfield-Kan integral
completion of \cite{bousfieldkan}, $\S$ I. It comes equipped with
a natural map $X_{\cdot}\rightarrow \Z_{\infty}X_{\cdot}$.

Following \cite{GilletSouleFiltrations}, $\S$3.1, we consider
$(\mathbf{T},\mathcal{O}_{\mathbf T})$ a ringed topos with
$\mathcal{O}_{\mathbf T}$ unitary and commutative. Then, for any
integer $N\geq 1$, the linear group of rank $N$ in $\mathbf{T}$,
$GL_N=GL_N(\mathcal{O}_{\mathbf T})$, is the sheaf associated to
the presheaf
$$U\mapsto GL_N(\Gamma(U,\mathcal{O}_{\mathbf T})).$$

Let $B_{\cdot}GL_N=B_{\cdot}GL_N(\mathcal{O}_{\mathbf T})$ be the
  classifying space of this sheaf of
groups. Observe that for every $N\geq
1$, there is a natural inclusion $B_{\cdot}GL_N\hookrightarrow
B_{\cdot}GL_{N+1}$. Consider the space $B_{\cdot}GL=\bigcup_N
B_{\cdot}GL_N$ and the following pointed spaces
\begin{eqnarray*}
\K_{\cdot} &=& \Z \times \Z_{\infty}B_{\cdot}GL, \\ \K_{\cdot}^N
&=& \Z \times \Z_{\infty}B_{\cdot}GL_N.
\end{eqnarray*}
Here, $\Z$ is the constant simplicial sheaf given by the constant
sheaf $\Z$, pointed by zero. For every $N\geq 1$, the direct sum
of matrices together with addition over $\Z$ gives a map
$$\K_{\cdot}^N\wedge \K_{\cdot}^N \rightarrow \K_{\cdot}.$$
These maps are compatible with the natural inclusions; thus
$\K_{\cdot}$ is equipped with an H-space structure (see
\cite{HuberWild}).

\subsection{K-theory }\label{ktheory}
Following  \cite{GilletSouleFiltrations}, for any space
$X_{\cdot}$, the \emph{stable $K$-theory}   is defined as
$$H^{-m}(X_{\cdot},\K_{\cdot})=[S^m\wedge X_{{\cdot}+},\K_{\cdot}],$$
and for every $N\geq 1$, the \emph{unstable $K$-theory} is defined
as
$$H^{-m}(X_{\cdot},\K_{\cdot}^N)=[S^m\wedge X_{{\cdot}+},\K_{\cdot}^N].$$
Since $\K_{\cdot}$ is an $H$-space, $H^{-m}(X_{\cdot},\K_{\cdot})$
are abelian groups for all $m$. However,
$H^{-m}(X_{\cdot},\K_{\cdot}^N)$ are abelian groups for all $m>0$
and in general only pointed sets for $m=0$.

\begin{df}
A space $X_{\cdot}$ is \emph{$K$-coherent}
if the natural maps
$$\lim_{\substack{\rightarrow \\ N}} H^{-m}(X_{\cdot},\K_{\cdot}^N)\rightarrow H^{-m}(X_{\cdot},\K_{\cdot})$$ and
$$\lim_{\substack{\rightarrow \\ N}} H^{m}(X_{\cdot},\pi_{n}\K_{\cdot}^N)\rightarrow H^{m}(X_{\cdot},\pi_{n}\K_{\cdot})$$
are isomorphisms for all $m,n\geq 0$.
\end{df}
(Here $H^{m}(X_{\cdot},\pi_{n}Y)$ are the singular cohomology
groups. See \cite{GilletSouleFiltrations}, $\S$1.2 for a
discussion in this language).

The Loday product induces a product structure on
$H^{-*}(X_{\cdot},\K_{\cdot})$ for every $K$-coherent space
$X_{\cdot}$.

\subsection{Comparison to Quillen's K-theory }
Let $(\mathbf{T},\mmO_{\mathbf T})$ be a locally ringed topos. For
every $U$ in $\mathbf{T}$, let $\mmP(U)$ be the category of
locally free $\mmO_{\mathbf{T}|U}$-sheaves of finite rank.

Let $B_{\cdot}Q\mmP$ be the simplicial sheaf obtained by the
Quillen construction applied to every $\mmP(U)$ (see
\cite{Quillen}). If $\Omega B_{\cdot}QP$ is the loop space of
$B_{\cdot}Q\mmP$, then, by the results of
\cite{GilletSouleFiltrations}, $\S$3.2.1 and \cite{Gillet},
Proposition 2.15, we obtain:
\begin{lema}\label{equiv}
In the homotopy category of simplicial sheaves, there is a natural map of
spaces
$$\Z \times \Z_{\infty}B_{\cdot}GL \rightarrow \Omega B_{\cdot}Q\mmP$$
which is a weak equivalence.
\end{lema}
Observe that this means that $\K_{\cdot}$ has yet another
$H$-space structure, given by the Waldhausen's pairing \cite{Waldhausen} $\Omega B_{\cdot}Q\mmP \wedge \Omega B_{\cdot}Q\mmP
\rightarrow \Omega B_{\cdot}Q\mmP$. As stated in \cite{GilletSouleFiltrations} \S 3.2.1, both structures agree.

It follows from the lemma that for any space $X_{\cdot}$ in
$\mathbf{sT}$, there is an isomorphism
$$H^{-m}(X_{\cdot},\K_{\cdot})\cong H^{-m}(X_{\cdot},\Omega B_{\cdot}Q\mmP).$$
Hence, the stable $K$-theory of a space can be computed  using the simplicial
sheaf $\Omega B_{\cdot}Q\mmP$ instead of the simplicial sheaf
$\K_{\cdot}$.

 Suppose that $\mathbf{T}$ is the category of sheaves over a
category of schemes $\mathbf{C}$. Let $\K_{\cdot}^{\sim}$ be a
fibrant resolution of $\Omega B_{\cdot}Q\mmP$. For every scheme
$X$ in $\mathbf{C}$, there is a natural map
\begin{equation}\label{quillen}
K_m(X)=\pi_m(\Omega B_{\cdot}Q\mmP(X))\rightarrow
\pi_m(\K_{\cdot}^{\sim}(X))\cong H^{-m}(X,\K_{\cdot}).
\end{equation}

The next theorem shows that many schemes are $K$-coherent and that
Quillen $K$-theory agrees with  stable $K$-theory.
\begin{theo}[\cite{GilletSouleFiltrations}, Proposition 5]\label{coherent}
Suppose that $X$ is a noetherian scheme of finite Krull dimension
$d$ and that $\mathbf{T}$ is either \begin{enumerate} \item
$\ZAR$, the big Zariski site of all noetherian schemes of finite
Krull dimension, \item $\ZAR(X)$, the big Zariski site of all
schemes of finite type over $X$, \item $\Zar(X)$, the small
Zariski site of $X$.
\end{enumerate}
Then, viewed as a $\mathbf{T}$-space, $X$ is $K$-coherent with
cohomological dimension at most $d$. Furthermore, the morphisms
$K_m(X)\rightarrow H^{-m}(X,\K_{\cdot})$ are isomorphisms for all
$m$.
\end{theo}

\begin{obs}
Let $\mmC$ be a small category
of schemes over $X$ that contains all open subschemes of its
objects. Consider the subsite $Z(X)$ of $\ZAR(X)$ obtained by
endowing $\mmC$ with the Zariski topology. Then, the statement of
the theorem will be true with $\mathbf{T}=T(Z(X))$.

For instance, if $X$ is a regular noetherian scheme of finite
Krull dimension, we could consider $Z(X)$ to be the site of all
regular schemes of finite type over $X$ as $Z(X)$. Another example
would be the site of all quasi-projective schemes of finite type
over a noetherian quasi-projective scheme of finite Krull
dimension.
\end{obs}

\subsection{$\mathbf{K}$-theory of spaces constructed from schemes}   Let
$\mathbf{C}=\ZAR$ and let $X_{\cdot}$ be a space constructed from
schemes. Then, in the Quillen context, one defines
$$K_{m}(X_{\cdot})=\pi_{m+1}(\holim\limits_{\substack{\longleftarrow \\ n}} B_{\cdot}
Q\mmP(X_n)).$$

For a description of the functor holim, see  \cite{bousfieldkan}, $\S$ XI, for
the case of simplicial sets or see \cite{Hirschhorn}, $\S 19$ for a general
treatment.

Observe that the construction of the map \eqref{quillen} can be
extended to spaces constructed from schemes. A space $X_{\cdot}$
is said to be \emph{degenerate} (above some simplicial degree) if
there exists an $N\geq 0$ such that $X_{\cdot}=\sk_N X_{\cdot}$
(where $\sk_N$ means the $N$-th skeleton of $X_{\cdot}$).

The next proposition is found in \cite{GilletSouleFiltrations},
$\S$3.2.3.
\begin{prop} Let $X_{\cdot}$ be a space constructed from schemes in $\ZAR$.
Then, the morphism \eqref{quillen} gives an isomorphism
$K_m(X_{\cdot})\cong H^{-m}(X_{\cdot},\K_{\cdot})$. Moreover, if
$X_{\cdot}$ is degenerate, then $X_{\cdot}$ is $K$-coherent.
\end{prop}

In particular, in the big Zariski site, since for every $N\geq 1$,
$B_{\cdot}GL_N$ is a simplicial scheme, we have
$K_{m}(B_{\cdot}GL_N)=H^{-m}(B_{\cdot}GL_N,\K_{\cdot})$. However,
$B_{\cdot}GL_N$ is not degenerate.

In a Zariski site over a base scheme $S$, the simplicial sheaf
$B_{\cdot}GL_N$ is the simplicial scheme given by the fibred
product $B_{\cdot}GL_{N/S}=B_{\cdot}GL_N \times_{\Z} S$.

\section{Characterization of maps from K-theory}
Our aim is to characterize functorial maps from $K$-theory. Since
stable $K$-theory is expressed as a representable functor, a first
approximation is obviously given by Yoneda's lemma. That is, given
a space $\F_{\cdot}$ and a map of spaces
$\K_{\cdot}\xrightarrow{\Phi} \F_{\cdot}$, the induced maps
$$H^{-m}(X_{\cdot},\K_{\cdot})\xrightarrow{\Phi_*} H^{-m}(X_{\cdot},\F_{\cdot}),\qquad \forall m\geq 0,$$
are determined by the image of $id \in  [\K_{\cdot},\K_{\cdot}]$
by the map $\Phi_*:H^0(\K_{\cdot},\K_{\cdot})\rightarrow
H^0(\K_{\cdot},\F_{\cdot})$. Indeed, if $g\in
H^{-m}(X_{\cdot},\K_{\cdot})=[S^m\wedge X_{\cdot},\K]$, there are
 induced morphisms
$$ [\K,\K]\xrightarrow{g^*} [S^m\wedge X_{\cdot},\K]\quad\textrm{and}\quad
 [\K,\F]\xrightarrow{g^*} [S^m\wedge X_{\cdot},\F].$$
Then, $g=g^*(id)$ and
$$\Phi_*(g)=\Phi_*g^*(id)=g^*\Phi_*(id). $$

We will see that, under some favorable conditions, the element $id$
can be changed by other universal elements at the level of the
simplicial scheme $B_{\cdot}GL_N$, for all $N\geq 1$.

\subsection{Compatible systems of maps and Yoneda lemma}
As in section \ref{K},  let
$(\mathbf{T},\mathcal{O}_{\mathbf{T}})$ be a ringed topos and let
$\F_{\cdot}$ be a fibrant space in $\mathbf{sT_*}$. A system of
maps $\Phi_M\in [\K_{\cdot}^M,\F_{\cdot}]$, $M\geq 1$, is said to
be \emph{compatible} if, for all $M'\geq M$, the diagram
$$
\xymatrix{\K_{\cdot}^{M} \ar[d] \ar[rd]^{\Phi_M} & \\
\K_{\cdot}^{M'} \ar[r]_{\Phi_{M'}} & \F_{\cdot} }
$$
is commutative in $\Ho(\mathbf{sT}_*)$. We associate to any map
$\Phi:\K_{\cdot}\rightarrow \F_{\cdot}$ in $\Ho(\mathbf{sT_*})$, a
compatible system of maps $\{\Phi_M\}_{M\geq 1}$,  given by the
composition of $\Phi$ with the natural map from $\K_{\cdot}^M$
into $\K_{\cdot}$.

Every compatible system of maps $\{\Phi_M\}_{M\geq 1}$ induces a natural
transformation of functors
$$\Phi(-): \lim_{\substack{\rightarrow\\
M}}[-,\K_{\cdot}^M] \rightarrow [-,\F_{\cdot}].$$

We state here a variant of Yoneda's lemma for maps induced by a compatible
system as above.

\begin{lema}\label{yoneda}
Let $\F_{\cdot}$ be a fibrant space in $\mathbf{sT_*}$. The map
$$\left\{\begin{array}{c}\textrm{compatible systems of maps } \\ \{\Phi_M\}_{M\geq
1},\  \Phi_M\in [\K_{\cdot}^M,\F_{\cdot} ] \end{array}  \right\}
\xrightarrow{\alpha}
\left\{\begin{array}{c} \textrm{natural transformation of functors }\\
\Phi(-): \lim_{\substack{\rightarrow\\
M}}[-,\K_{\cdot}^M] \rightarrow [-,\F_{\cdot}] \end{array}
\right\}$$ sending every compatible system of maps to its induced
natural transformation, is a bijection.
\end{lema}
\begin{proof}
We prove the result by giving the explicit inverse arrow $\beta$
of $\alpha$.

So, let
$$\Phi(-): \lim_{\substack{\rightarrow\\ M}}[-,\K_{\cdot}^M] \rightarrow
[-,\F_{\cdot}]$$ be a natural transformation of functors. For
every $N\geq 1$, let $e_N\in \lim_{\substack{\rightarrow \\
M}}[\K_{\cdot}^N,\K_{\cdot}^M]$ be the image of $id \in
[\K_{\cdot}^N,\K_{\cdot}^N]$ under the natural morphism
$ [\K_{\cdot}^N,\K_{\cdot}^N] \xrightarrow{\sigma_N} \lim_{\substack{\rightarrow \\
M}}[\K_{\cdot}^N,\K_{\cdot}^M].$
 We define $\Phi_N=\beta(\Phi)_N\in [\K_{\cdot}^N,\F_{\cdot}]$ to be the image of $e_N$ by
$\Phi$, $$\Phi_N := \Phi(\K_{\cdot}^N)(e_N). $$ In general, for
every $N'\geq N\geq 1$, consider the map $e_{N,N'}\in
[\K_{\cdot}^N,\K_{\cdot}^{N'}]$ induced by the natural inclusion
$B_{\cdot}GL_N \hookrightarrow B_{\cdot}GL_{N'}.$ Observe that the
image of $e_{N,N'}$ under the   map
$$ [\K_{\cdot}^N,\K_{\cdot}^{N'}] \xrightarrow{\sigma_{N'}} \lim_{\substack{\rightarrow \\
M}}[\K_{\cdot}^N,\K_{\cdot}^M],$$ is exactly $e_N$. Moreover,  by
hypothesis, there is a commutative diagram
$$
\xymatrix{
\lim_{\substack{\rightarrow\\
M}}[\K_{\cdot}^{N'},\K_{\cdot}^M] \ar[d]_{e_{N,N'}^*} \ar[rr]^(.6){\Phi(\K_{\cdot}^{N'})} && [\K_{\cdot}^{N'},\F_{\cdot}] \ar[d]^{e_{N,N'}^*}\\
\lim_{\substack{\rightarrow\\
M}}[\K_{\cdot}^{N},\K_{\cdot}^M]
\ar[rr]_(.6){\Phi(\K_{\cdot}^N)}&& [\K_{\cdot}^{N},\F_{\cdot}] }
$$ which gives the compatibility of the system $\{\Phi_N\}_{N\geq
1}$. Therefore, the map $\beta$ is defined.

Now let $X_{\cdot}$ be any space in $\mathbf{sT_*}$. In order to
prove that $\beta$ is a right inverse of $\alpha$, we have to see
that $\Phi(X_{\cdot})$ is the map induced by the just constructed
system $\{\Phi_M\}_{M\geq 1}$.

Let $f\in \lim_{\substack{\rightarrow \\ M}}
[X_{\cdot},\K_{\cdot}^M]$. Then,  there exists an integer $N\geq
1$ and a map $g\in [X_{\cdot},\K_{\cdot}^N]$, such that
$\sigma_N(g)=f$. By the commutative diagram
$$\xymatrix{
[\K_{\cdot}^N,\K_{\cdot}^N] \ar[r]^(.45){\sigma_N} \ar[d]_{g^*} & \lim_{\substack{\rightarrow \\
M}} [\K_{\cdot}^N,\K_{\cdot}^M] \ar[d]^{g^*} \\
[X_{\cdot},\K_{\cdot}^N]
\ar[r]_(.45){\sigma_N}   & \lim_{\substack{\rightarrow \\
M}}[X_{\cdot},\K_{\cdot}^M], } $$ we see that in fact,
$f=\sigma_N(g)=g^*(e_N).$ Using the fact that $\Phi$ is a natural
transformation, the diagram
\begin{equation}\label{diagram3}\begin{array}{c}\xymatrix{ \lim_{\substack{\rightarrow \\ M}}
[\K_{\cdot}^{N},\K_{\cdot}^{M}] \ar[r]^(.6){\Phi} \ar[d]_{g^*} &
[\K_{\cdot}^{N},\F_{\cdot}] \ar[d]^{g^*}
\\ \lim_{\substack{\rightarrow \\ M}}[X_{\cdot},\K_{\cdot}^{M}]
\ar[r]_(.6){\Phi} & [X_{\cdot},\F_{\cdot}]
}\end{array}\end{equation}  is commutative.  Hence, we obtain
$$\Phi(f)=\Phi(g^*(e_N))=g^*\Phi(e_N)=\beta(\Phi)_N\circ g =\alpha\beta(\Phi)(f),$$
as desired.

It remains to check that $\beta$ is a left inverse of $\alpha$.
Let $\{\Phi_N\}_{N\geq 1}$ be a compatible system of maps, let
$\Phi$ be the associated transformation of functors obtained by
$\alpha$ and let $\{\beta(\Phi)_N\}_{N\geq 1}$ be the system
$\beta(\Phi)$. From the commutative diagram
$$ \xymatrix{
 [\K_{\cdot}^N,\K_{\cdot}^N] \ar[d]_{\sigma_N} \ar[r]^{(\Phi_N)_*}   & [\K_{\cdot}^N,\F_{\cdot}] \\
 \lim_{\substack{\rightarrow \\ M}} [\K_{\cdot}^{N},\K_{\cdot}^{M}] \ar[ur]_{\Phi}
}$$ we deduce that
$$\Phi_{N}=(\Phi_N)_*(id)=\Phi\sigma_N(id)=\Phi(e_N)=\beta(\Phi)_N . $$
Therefore, $\beta$ is the inverse of $\alpha$ and thus $\alpha$ is
a bijection.
\end{proof}

\begin{obs}
The last lemma is not specific to our category and to our
compatible system of maps. It could be directly generalized to any
suitable category.
\end{obs}

\subsection{Weakly additive systems of maps}\label{weak} We start by defining
the class of  weakly additive systems of maps. It is for this class of maps that we
will state our results on the comparison of the induced maps. It will be shown below that many usual maps are
weakly additive.

Let $\K_{\cdot}^N,\K_{\cdot}$  be as in the previous section. Let
$\pr_1$ and $\pr_2$ be the projections onto the first and second
component respectively
\begin{eqnarray*}
\pr\nolimits_1 & : & \Z\times \Z_{\infty}B_{\cdot}GL_M \rightarrow
\Z,
\\ \pr\nolimits_2 & : & \Z\times \Z_{\infty}B_{\cdot}GL_M
\rightarrow \Z_{\infty}B_{\cdot}GL_M, \end{eqnarray*} and let
$j_1,j_2$ denote the inclusions obtained using the respective base
points
\begin{eqnarray*}
\Z & \xrightarrow{j_1} & \Z\times \Z_{\infty}B_{\cdot}GL_M,
\\ \Z_{\infty}B_{\cdot}GL_M & \xrightarrow{j_2} & \Z\times
\Z_{\infty}B_{\cdot}GL_M.
\end{eqnarray*}
Denote by $\pi_i=j_i\circ \pr_i\in [\K_{\cdot}^M,\K_{\cdot}^M]$,
$i=1,2$, the compositions
$$
\begin{array}{cccccc}
\pi_1: &  \Z\times \Z_{\infty}B_{\cdot}GL_M & \xrightarrow{\pr_1}
& \Z & \xrightarrow{j_1} & \Z\times \Z_{\infty}B_{\cdot}GL_M,
\\ \pi_2: & \Z\times \Z_{\infty}B_{\cdot}GL_M & \xrightarrow{\pr_2} &
\Z_{\infty}B_{\cdot}GL_M & \xrightarrow{j_2} & \Z\times
\Z_{\infty}B_{\cdot}GL_M.
\end{array}
$$

For every space $\F_{\cdot}$ in $\mathbf{sT_*}$, there are induced
maps
\begin{eqnarray*}
[\K_{\cdot}^M,\F_{\cdot}]& \xrightarrow{\pi_i^*} &
[\K_{\cdot}^M,\F_{\cdot}],\qquad i=1,2.
\end{eqnarray*}
If $\Phi_M\in [\K_{\cdot}^M,\F_{\cdot}]$ we define the maps
$$\Phi_M^i:= \pi_i^*(\Phi_M)\in [\K_{\cdot}^M,\F_{\cdot}],\qquad i=1,2.$$

\begin{obs}\label{constant} Consider a compatible system of
maps $\Phi_M\in [\K_{\cdot}^M,\F_{\cdot}]$, $M\geq 1$ and assume
that $\F_{\cdot}$ is fibrant (this is no loss of generality). By
hypothesis,  the diagrams
$$\begin{array}{ccc}\xymatrix{& \K_{\cdot}^M \ar[d] \\ \Z \ar[r]_{j_1} \ar[ru]^{j_1} & \K_{\cdot}^{M'}} & \qquad
\xymatrix{\\ \textrm{and}}\qquad & \xymatrix{\K_{\cdot}^{M} \ar[d] \ar[rd]^{\Phi_M} & \\
\K_{\cdot}^{M'} \ar[r]_{\Phi_{M'}} & \F_{\cdot} }
\end{array}$$  are commutative in $\Ho(\mathbf{sT}_*)$, for every $M'\geq M$. Therefore,
the homotopy class of the map $\Phi_M \circ j_1$ in
$[\Z,\F_{\cdot}]$ does not depend on $M$.
\end{obs}

\begin{df}\label{weakly} Let $\F_{\cdot}$ be any space in $\mathbf{sT_*}$.
By an operation $\bullet=\{\bullet_M\}_{M\geq 1}$ on $\{[\K^M_{\cdot},\F_{\cdot}]\}_{M\geq 1}$ we mean a collection of maps
$$\bullet_M :  [\K_{\cdot}^M,\F_{\cdot}] \times [\K_{\cdot}^M,\F_{\cdot}] \rightarrow [\K_{\cdot}^M,\F_{\cdot}],\quad M\geq 1.$$
Note that we do not require
compatibility for different indices $M$.

Let $\{\Phi_M\}_{M\geq 1}$, with $\Phi_M\in
[\K_{\cdot}^M,\F_{\cdot}]$, be a compatible system of maps.  The system
$\{\Phi_M\}_{M\geq 1}$ is called \emph{weakly additive with
respect to an operation $\bullet=\{\bullet_M\}_{M\geq 1}$}, if
for every $M$, there is an equality
$$\Phi_M= \Phi^1_M \bullet_M \Phi_M^2. $$

Any map  $\bullet:[\K_{\cdot},\F_{\cdot}] \wedge [\K_{\cdot},\F_{\cdot}] \rightarrow [\K_{\cdot},\F_{\cdot}]$ induces an operation on $\{[\K^M_{\cdot},\F_{\cdot}]\}_{M\geq 1}$
by composition with the natural maps $\K^M_{\cdot}\rightarrow \K_{\cdot}$. Then, a map $\Phi\in [\K_{\cdot},\F_{\cdot}]$
is called \emph{weakly
additive} with respect to $\bullet$ if the induced compatible system of maps is weakly
additive with respect to the induced operation $\{\bullet_M\}_{M\geq 1}$.
\end{df}

\begin{ex}[$\Phi$ is trivial over $\Z$]\label{trivial2} Let $\{\Phi_M\}_{M\geq 1}$ be a compatible
system of maps with $\Phi_M\in [\K_{\cdot}^M,\F_{\cdot}]$. Assume
that for all $M\geq 1$, $\Phi_M^1=*$, i.e. the constant map to the
base point of $\F_{\cdot}$. Then, $\Phi_M=\Phi^2_M$ for all $M$.
Therefore, if we take the operation $\bullet_M=\pr\nolimits_2$, the system
$\{\Phi_M\}_{M\geq 1}$ is weakly additive with respect to $\bullet=\{\bullet_M\}_{M\geq 1}$.
\end{ex}

\begin{ex}[$\F_{\cdot}$ is an H-space]\label{hspace} Assume that $\F_{\cdot}$ is an $H$-space. Then, one can take the operation $\bullet_M$
to be the sum in $[\K_{\cdot}^M,\F_{\cdot}]$.
Then for a compatible system of maps $\{\Phi_M\}_{M\geq 1}$  the condition of weakly additivity  means that the maps $\Phi_M$ behave additively over the two components of
$\K_{\cdot}^M$. Actually, the definition of weakly additive
systems of maps was motivated by this example.

The lambda operations on higher algebraic $K$-theory, defined by
Gillet and Soul\'e in \cite{GilletSouleFiltrations}, are an
example of this type of weakly additive systems of maps.
\end{ex}

\begin{obs}\label{trivialhspace}
If $\F_{\cdot}$ is an $H$-space and $\Phi=\{\Phi_M\}_{M\geq 1}$ is a
compatible system of maps such that  $\Phi_M^1=*$, then the system $\Phi$
is weakly additive with respect to the $H$-sum of
$\F_{\cdot}$ (previous example) and also with respect to the operation $\bullet_M=\pr_2$ (by example \ref{trivial2}).
\end{obs}

\subsection{Classifying elements} We now introduce some classifying elements.
For every $N\geq 1$, let
$$ \sigma_N: [\Z_{\infty}B_{\cdot}GL_N,\K_{\cdot}^N] \rightarrow  \lim_{\substack{\rightarrow \\
M}}[\Z_{\infty}B_{\cdot}GL_N,\K_{\cdot}^M]$$ be the natural
morphism. Recall that we defined $j_2\in
[\Z_{\infty}B_{\cdot}GL_N,\K_{\cdot}^N]$ to be  the map induced by
the natural inclusion. Then, we define
$$i_N'=\sigma_N(j_2) \in   \lim_{\substack{\rightarrow \\
M}}[\Z_{\infty}B_{\cdot}GL_N,\K_{\cdot}^M].$$
Let $r\geq 0$ and $u'_{r}\in \lim_{\substack{\rightarrow \\
M}}[\Z_{\infty}B_{\cdot}GL_N,\K_{\cdot}^M]$ be the image by $\sigma_N$ of the homotopy class of the constant map
\begin{eqnarray*}
\Z_{\infty}B_{\cdot}GL_N &\rightarrow & \Z\times \Z_{\infty}B_{\cdot}GL_N \\
x & \mapsto & (r,*).
\end{eqnarray*}
Finally, consider the natural map $B_{\cdot}GL_N\rightarrow
\Z_{\infty}B_{\cdot}GL_N$.
The images of $i_N'$ and $u_r'$  under the induced map $$\lim_{\substack{\rightarrow \\
M}}[\Z_{\infty}B_{\cdot}GL_N,\K_{\cdot}^M]\rightarrow \lim_{\substack{\rightarrow \\
M}}[B_{\cdot}GL_N,\K_{\cdot}^M],$$
are denoted by $$i_N,u_r\in \lim_{\substack{\rightarrow \\
M}}[B_{\cdot}GL_N,\K_{\cdot}^M]$$ respectively.

\begin{prop}\label{uniqueness1}
Let $\F_{\cdot}$ be a space in $\mathbf{sT_*}$ and let
$\{\Phi_M\}_{M\geq 1}$, $\{\Phi'_M\}_{M\geq 1}$ be two weakly
additive systems of maps with respect to the same operation. Then,
the induced maps
$$\Phi,\Phi': \lim_{\substack{\rightarrow\\ M}}[-,\K_{\cdot}^M]
\rightarrow [-,\F_{\cdot}]$$ agree for all spaces, if and only if,
in $[\Z_{\infty}B_{\cdot}GL_N,\F_{\cdot}]$ it holds
\begin{eqnarray}
\Phi(i_N')&=& \Phi'(i_N'), \quad \textrm{for all } N\geq 1, \label{phi2}\\
\Phi(u'_r)&=& \Phi'(u'_r),\quad \textrm{for all } r\in \Z,\ N\geq 1.
\label{phi1}
\end{eqnarray}
\end{prop}

\begin{proof} One implication is obvious. By lemma \ref{yoneda}, it is enough to see that for all $N$,
$\Phi_N=\Phi_N'$. By hypothesis, there is an operation $\bullet_N$
on $[\K_{\cdot}^N,\F_{\cdot}]$ such that
\begin{eqnarray*}
\Phi_N& =& \Phi_N^1 \bullet_N \Phi_N^2, \\
\Phi_N'& =& \Phi_N^{'1} \bullet_N \Phi_N^{'2}.
\end{eqnarray*}
Therefore, it is enough to see that $$ \Phi_N^1 = \Phi_N^{'1},\quad \textrm{and}
\quad \Phi_N^{2}=\Phi_N^{'2}.$$ The first equality follows from hypothesis
\eqref{phi1}. For the second equality, observe that by definition,
$\Phi(i_N')=\Phi_N\circ j_2$. Therefore, by equality \eqref{phi2},
$$\Phi^2_N=\Phi_N\circ j_2\circ \pr\nolimits_2 = \Phi(i_N')\circ\pr\nolimits_2
=\Phi'(i_N')\circ\pr\nolimits_2=\Phi^{'2}_N.$$
\end{proof}

\begin{cor}\label{uniqueness2}
Let $\F_{\cdot}$ be an $H$-space in $\mathbf{sT_*}$ and
$\{\Phi_M\}_{M\geq 1}$, $\{\Phi'_M\}_{M\geq 1}$ be two weakly
additive systems of maps with respect to
 the same operation.
Then, the induced maps
$$\Phi(-),\Phi'(-): \lim_{\substack{\rightarrow\\ M}}[-,\K_{\cdot}^M]
\rightarrow [-,\F_{\cdot}]$$ agree if and only if, in
$[B_{\cdot}GL_N,\F_{\cdot}]$ it  holds
\begin{eqnarray}
\Phi(i_N)&=& \Phi'(i_N), \quad \textrm{for all } N\geq 1, \label{phi22}\\
\Phi(u_r)&=& \Phi'(u_r),\quad \textrm{for all } r\in \Z,\ N\geq 1.
\label{phi11}
\end{eqnarray}
\end{cor}
\begin{proof}
It follows from the fact that the natural map
$$ [\Z_{\infty}B_{\cdot}GL_N,\F_{\cdot}]\rightarrow [B_{\cdot}GL_N,\F_{\cdot}]$$ is an isomorphism if $\F_{\cdot}$ is an $H$-space,
 and under this isomorphism, the elements $u_r$ and $i_N$
correspond to $u_r'$ and $i_N'$ respectively.
\end{proof}

Let $j_N\in H^0(B_{\cdot}GL_{N},\K_{\cdot})=K_0(B_{\cdot}GL_{N})$
be the image of $i_N$ under the morphism
$$\lim_{\substack{\rightarrow \\ M}} H^{0}(B_{\cdot}GL_{N},\K_{\cdot}^M) \rightarrow H^0(B_{\cdot}GL_{N},\K_{\cdot}).$$
By abuse of notation, we denote by  $u_r$ the image of $u_r\in \lim_{\substack{\rightarrow
\\ M}} H^{0}(B_{\cdot}GL_{N},\K_{\cdot}^M)$ in $H^0(B_{\cdot}GL_{N},\K_{\cdot})$.

\begin{cor}\label{uniqueness3}
Let $\F_{\cdot}$ be an $H$-space in $\mathbf{sT_*}$. Let
$\Phi,\Phi':\K_{\cdot}\rightarrow \F_{\cdot}$ be two H-space maps.
Then, $\Phi$ and $\Phi'$ are weakly additive. Moreover,
$$\Phi,\Phi':[X_{\cdot},\K_{\cdot}]\rightarrow [X_{\cdot},\F_{\cdot}],$$ agree for all $K$-coherent spaces $X_{\cdot}$, if they agree at $j_N,u_r \in K_0(B_{\cdot}GL_N)$  for all $N\geq 1$ and all $r\in \Z$.
\end{cor}
\begin{proof} The maps $\Phi$ and $\Phi'$ are weakly additive with respect to the $H$-sum of $\F_{\cdot}$, due to example
\ref{hspace}. Therefore, by corollary \ref{uniqueness2}, the maps
$$\Phi,\Phi':[X_{\cdot},\K_{\cdot}]\cong \lim_{\substack{\rightarrow\\ N}}[X_{\cdot},\K_{\cdot}^N]
\rightarrow [X_{\cdot},\F_{\cdot}]$$ agree for all $K$-coherent
spaces $X_{\cdot}$, if and only if $\Phi(i_N)=\Phi'(i_N)$ for
$N\geq 1$ and $\Phi(u_r)=\Phi'(u_r)$ for all $r$. Since by
construction $\Phi(i_N)=\Phi(j_N)$ (and the same for $\Phi'$), the
corollary is proved.
\end{proof}

\subsection{Application to the Zariski sites}\label{zariski} Let $S$ be a  finite dimensional noetherian scheme.
Fix $\mathbf{C}$ a Zariski subsite of $\ZAR(S)$ containing all
open subschemes of its objects and the components of the simplicial scheme
$B_{\cdot}GL_{N/S}$. Let $\mathbf{T}=T(\mathbf{C})$.

A direct consequence of corollary \ref{uniqueness2} is the
following theorem.
\begin{theo}\label{uni7} Let $\F_{\cdot}$ be an $H$-space in $\mathbf{sT}_*$  and $\{\Phi_M\}_{M\geq 1}$,
$\{\Phi'_M\}_{M\geq 1}$ be two weakly additive systems of maps
with respect to the same   operation. Then,  the induced maps
$$\Phi,\Phi': K_m(X_{\cdot})\cong \lim_{\substack{\rightarrow\\ M}} H^{-m}(X_{\cdot},\K_{\cdot}^M)
\rightarrow H^{-m}(X_{\cdot},\F_{\cdot})$$ agree for all $m\geq 0$
and all $K$-coherent spaces $X_{\cdot}$, if and only if in
$H^0(B_{\cdot}GL_{N/S},\F_{\cdot})$ it holds
\begin{eqnarray*}
\Phi(i_N)&=& \Phi'(i_N), \quad \textrm{for all } N\geq 1, \\
\Phi(u_r)&=& \Phi'(u_r),\quad \textrm{for all } r\in \Z,\ N\geq 1.
\end{eqnarray*}
\end{theo}
\begin{flushright}
$\square$
\end{flushright}

Finally, the next corollary is corollary \ref{uniqueness3} applied
to the Zariski subsite $\mathbf{C}\subset \ZAR(S)$.

\begin{cor}\label{uni} Let $\F_{\cdot}$ be an $H$-space in $\mathbf{sT}_*$.
Let $\chi_1,\chi_2\in [\K_{\cdot},\F_{\cdot}]$ be two $H$-space
maps in $\Ho(\mathbf{sT}_*)$. Then, the induced maps
$$\chi_1,\chi_2: K_m(X_{\cdot}) \rightarrow H^{-m}(X_{\cdot},\F_{\cdot}) $$
agree for all degenerate simplicial schemes in $\mathbf{C}$,  if for all $N\geq
1$ and $r\in \Z$, they agree at $j_N,u_r \in K_0(B_{\cdot}GL_{N/S})$.
\end{cor}
\begin{flushright}
$\square$
\end{flushright}

The next theorem shows that in fact a weaker condition is needed in order to
obtain the uniqueness of maps.

For any scheme $X$ and any simplicial sheaf $\F_{\cdot}$, let
$\F_{X\cdot}$ denote the restriction of $\F_{\cdot}$ to the small
Zariski site of $X$. In the next theorem, we write
$[\cdot,\cdot]_{\mathbf{C}}$ for the maps in
$\Ho(\mathbf{sT(C)}_*)$, for any site $\mathbf{C}$. If $X$ is a
scheme in $\mathbf{C}$, let $\mathbf{C}(X)$
 be the subsite of $\ZAR(X)$ whose objects are in $\mathbf{C}$.

\begin{theo}\label{uni2}
Let $\F_{\cdot}$ be a pseudo-flasque sheaf on $\mathbf{sT}_*$,
which is an $H$-space. Assume that
\begin{enumerate*}[$\blacktriangleright$]
\item For every scheme $X$ in $\mathbf{C}$, there are two
$H$-space maps $$\chi_i(X)\in
[\K_{X\cdot},\F_{X\cdot}]_{\Zar(X)},\qquad i=1,2.$$ \item For any map
$X\rightarrow Y$ in $\mathbf{C}$ and for $i=1,2$ there is a commutative diagram
\begin{equation}\label{Fcommut} \xymatrix{
\K_{\cdot}(Y) \ar[rr]^{\chi_i(Y)}\ar[d] && \F_{\cdot}(Y) \ar[d] \\
\K_{\cdot}(X) \ar[rr]^{\chi_i(X)}  && \F_{\cdot}(X) }
\end{equation} in $\Ho(\mathbf{SSets}_*)$. \end{enumerate*} Then,
the maps
$$\chi_1,\chi_2: K_m(X) \rightarrow H^{-m}(X,\F_{\cdot}) $$
agree for all  schemes in $\mathbf{C}$ if they agree at $j_N,u_r
\in K_0(B_{\cdot}GL_{N/S})$, for all $N\geq 1$ and $r\in \Z$ (see
remark below).
\end{theo}

\begin{obs}\label{pseudoflasque}
The condition of $\F_{\cdot}$ being pseudo-flasque means that
for any space $X_{\cdot}$ constructed from schemes,
$$H^{-m}(X_{\cdot},\F_{\cdot})\cong \pi_m(\F_{\cdot}(X_{\cdot})):=
\pi_m(\holim\limits_{\substack{\longleftarrow \\ n}} \F_{\cdot}(X_n)). $$
Now, the commutative diagram \eqref{Fcommut}, implies that, for any
such space, there are induced morphisms
$$\K_{\cdot}(X_{\cdot}) \xrightarrow{\chi_i(X_{\cdot})} \F_{\cdot}(X_{\cdot}),\qquad i=1,2. $$
Hence, the maps $\chi_1$ and $\chi_2$ are defined for
$X_{\cdot}=B_{\cdot}GL_{N/S}$.
\end{obs}

\begin{proof}
For every fixed scheme $X$, it follows from theorem \ref{uni}
  that $\chi_1(X)=\chi_2(X)$ if they agree at $j_N,u_r\in
[B_{\cdot}GL_{N|X},\K_{X\cdot}]_{\Zar(X)}$. Observe now that by
the remarks following proposition \ref{coherent},
$$[B_{\cdot}GL_{N|X},\K_{X\cdot}]_{\Zar(X)}\cong
[B_{\cdot}GL_{N|X},\K_{\cdot}]_{\mathbf{C}(X)}=K_0(B_{\cdot}GL_{N|X}),$$
and for $i=1,2$ there is a commutative diagram
\begin{equation}\label{pseudoflasque2}
\xymatrix{ K_0(B_{\cdot}GL_{N|S}) \ar[rr]^{\chi_i} \ar[d]& &[B_{\cdot}GL_{N|S},\F_{\cdot}]_{\mathbf{C}} \ar[d]\\
K_0(B_{\cdot}GL_{N|X})\ar[rr]_{\chi_i} &&
[B_{\cdot}GL_{N|X},\F_{\cdot}]_{\mathbf{C}(X)} }
\end{equation}
Then, the statement follows from the fact that $j_N$ and $u_r$ in
$K_0(B_{\cdot}GL_{N|X}) $ are the image under the vertical map of
$j_N$ and $u_r$ in $K_0(B_{\cdot}GL_{N|S})$.
\end{proof}

\section{Morphisms between K-groups }

\subsection{Lambda and Adams operations}
In this section we focus on the case where
$\F_{\cdot}=\K_{\cdot}$. Then, the main application of theorems
\ref{uni} and \ref{uni2} is to the Adams operations and to the
lambda operations on higher algebraic $K$-theory.

The
Grothendieck group of a scheme $X$, has a $\lambda$-ring structure
given by $\lambda^k(E)=\bigwedge^k E$, for any vector bundle $E$
over $X$.
 In the literature there
are several definitions of the extension of the Adams operations
of $K_0(X)$ to the higher algebraic K-groups. Our aim in this
section is to give a criterion for their comparison.

Soul{\'e}, in \cite{Soule}, gives a $\lambda$-ring structure to the higher
algebraic $K$-groups of any noetherian regular scheme of finite Krull dimension.
Gillet and Soul{\'e} then generalize this result in \cite{GilletSouleFiltrations},
defining lambda operations for all $K$-coherent spaces in any locally ringed
topos. We briefly recall this construction here.

Let  $\mathbf{R}_{\Z}(GL_N)$ be  the Grothendieck group of
representations of the general linear group scheme $GL_{N/\Z}$.
The properties of $\mathbf{R}_{\Z}(GL_N)$ that concern us are:
\begin{enumerate}[(1)]
\item $\mathbf{R}_{\Z}(GL_N)$ has a $\lambda$-ring structure.  \item For any
locally ringed topos, there is a ring morphism
$$\varphi:\mathbf{R}_{\Z}(GL_N)\rightarrow H^0(B_{\cdot}GL_N,\K_{\cdot}).$$
\end{enumerate}

The operations $\lambda^k_{GS},\Psi^k_{GS}$ are constructed by
transferring the lambda and Adams operations of
$\mathbf{R}_{\Z}(GL_N)$ to the $K$-theory of $B_{\cdot}GL_N$.
Namely, consider the representation $id_N-N$ and the maps
$$\varphi(\Psi^k(id_N-N)),\ \varphi(\lambda^k(id_N-N)) : B_{\cdot}GL_N \rightarrow \K_{\cdot} $$
 in the homotopy
category of simplicial sheaves.

Consider the only $\lambda$-ring structure on $\Z$ with trivial
involution. Then, adding the previous maps with the lambda or
Adams operations on the $\Z$-component, we obtain compatible
systems of maps
$$\Psi_{GS}^k,\ \lambda_{GS}^k: \K_{\cdot}^N\rightarrow \K_{\cdot},\qquad N\geq 1. $$
Observe that, by example \ref{hspace}, both systems are weakly
additive with respect to the $H$-sum of $\K_{\cdot}$.

In particular, for any noetherian scheme $X$ of finite Krull dimension, there
are induced Adams operations on the higher $K$-groups
$$\Psi^k_{GS}: K_m(X) \rightarrow K_m(X). $$
Gillet and Soul\'e checked that these maps satisfy the identities
of a special lambda ring.

\subsection{Vector bundles over a simplicial scheme }
Let $X_{\cdot}$ be a simplicial scheme, with face maps denoted by
$d_i$ and degeneracy maps by $s_i$. A \emph{vector bundle}
$E_{\cdot}$ over $X_{\cdot}$ consists of a collection of vector
bundles $E_n\rightarrow X_n$, $n\geq 0$, together with
isomorphisms $d_i^*E_{n}\cong E_{n+1}$ and $s_i^*E_{n+1}\cong E_n$
for all face and degeneracy maps. Moreover, these isomorphisms
should satisfy the simplicial identities. By a \emph{morphism of
vector bundles} we mean a collection of morphisms at each level,
compatible with these isomorphisms. An \emph{exact sequence} of
vector bundles is an exact sequence at every level.

Let $\Vect(X_{\cdot})$ be the exact category of vector bundles
over $X_{\cdot}$ and consider the algebraic $K$-groups of
$\Vect(X_{\cdot})$, $K_m(\Vect(X_{\cdot}))$. These can be computed as the homotopy groups of
the simplicial set $S_{\cdot}(\Vect(X_{\cdot}))$ given by the Waldhausen
construction.

For every simplicial scheme $X_{\cdot}$ and every $n\geq 0$, there
is a natural simplicial map
$$S_{\cdot}(\Vect(X_{\cdot})) \rightarrow  S_{\cdot}(\Vect(X_n)).$$
By the definition of  vector bundles over  simplicial schemes, it induces a
simplicial map
$$S_{\cdot}(\Vect(X_{\cdot})) \rightarrow \holim_{\substack{\longleftarrow \\ n}} S_{\cdot}
(\Vect(X_n)),$$ which induces a
morphism
$$K_m(\Vect(X_{\cdot})) \xrightarrow{\psi} K_m(X_{\cdot}),\qquad m\geq 0.$$

At the zero level, $K_0(\Vect(X_{\cdot}))$ is the Grothendieck
group of the category of vector bundles over $X_{\cdot}$, and
hence it has a $\lambda$-ring structure.

In the particular situation where $X_{\cdot}$ is a simplicial
object in $\ZAR(S)$, with $S$ a  finite dimensional noetherian
scheme, the above morphism can be described as follows (see
\cite{GilletSouleFiltrations}).

Let $N_{\cdot}U$ denote the nerve of a covering $U$ of $X$. Let
$E^M_{\cdot}$ denote the universal vector bundle over
$B_{\cdot}GL_{M/S}$ and let $E_{\cdot}$ be a rank $N$ vector
bundle over $X_{\cdot}$. Then, there exists a hypercovering
$p:N_{\cdot}U\rightarrow X_{\cdot}$  and a classifying map $\chi:
N_{\cdot}U\rightarrow B_{\cdot}GL_{M/S}$, for $M\geq N$, such that
$p^*(E_{\cdot})=\chi^*(E_{\cdot}^M)$. The induced map
$$\chi: N_{\cdot}U\rightarrow \{N\}\times \Z_{\infty}B_{\cdot}GL_M \rightarrow \Z\times
\Z_{\infty}B_{\cdot}GL_M\rightarrow \K_{\cdot}$$ in $\ZAR(S)$,
defines an element $\chi$ in
$H^0(N_{\cdot}U,\K_{\cdot})=H^0(X_{\cdot},\K_{\cdot})=K_0(X_{\cdot})$,
which is $\psi([E_{\cdot}])$. This description also shows that the
morphism factorizes through the limit
$$\psi:K_0(\Vect(X_{\cdot}))\rightarrow \lim_{\substack{\rightarrow \\ M} } H^0(X_{\cdot},\K_{\cdot}^M)
\rightarrow H^0(X_{\cdot},\K_{\cdot}).$$ When
$X_{\cdot}=B_{\cdot}GL_{N/S}$, we obtain that
\begin{eqnarray*}
\psi(E^N_{\cdot}-N)&=& j_N \in K_0(B_{\cdot}GL_{N/S}), \\ \psi(E^N_{\cdot}-N) & = & i_N \in \lim_{\substack{\rightarrow \\
M} } H^0(B_{\cdot}GL_{N/S},\K_{\cdot}^M).
\end{eqnarray*}
Here $N$ is the trivial bundle of rank $N$. Clearly, the trivial
bundle of rank $r\geq 0$ in $B_{\cdot}GL_{N/S}$, is mapped to
$u_r$.

Consider the $\lambda$-ring structure in
$K_0(\Vect(B_{\cdot}GL_{N/S}))$ and denote by $\Psi^k$ the
corresponding Adams operations. Gillet and Soul{\'e} proved in
\cite{GilletSouleFiltrations} section 5 that there are equalities
\begin{eqnarray*}
\varphi(\Psi^k(id_N-N))&=& \psi(\Psi^k(E_{\cdot}^N-N)),\\
\varphi(\lambda^k(id_N-N))&=&\psi(\lambda^k(E_{\cdot}^N-N)).
\end{eqnarray*}
 Moreover, one can
easily check that $\varphi(\Psi^k(id_N-N))=\Psi^k_{GS}(i_N)$, and
$\varphi(\lambda^k(id_N-N))=\lambda^k_{GS}(i_N)$. Therefore,
$$\Psi^k_{GS}(i_N)=\psi(\Psi^k(E_{\cdot}^N-N)),\quad\textrm{and}\quad
\lambda^k_{GS}(i_N)=\psi(\lambda^k(E_{\cdot}^N-N)).$$ Also, it
holds by definition that
$$
\Psi^k_{GS}(u_r)=\psi(\Psi^k(u_r))=u_{\psi^k(r)},\quad\textrm{and}\quad\lambda^k_{GS}(u_r)=\psi(\lambda^k(u_r))=
u_{\lambda^k(r)}.
$$

\subsection{Uniqueness theorems} Let $S$ be a finite dimensional noetherian scheme.
Fix $\mathbf{C}$ a Zariski subsite of $\ZAR(S)$ as in section
\ref{zariski}. The following theorems are a consequence of
theorems  \ref{uni} and \ref{uni2}  applied to the present
situation.

\begin{theo}[Lambda operations]\label{uni4} Let $\{\rho_N:\K_{\cdot}^N\rightarrow \K_{\cdot}\}_{N\geq 1}$ be a weakly additive
system of maps with respect to the $H$-sum of $\K_{\cdot}$. Let
$\rho$ be the induced morphism
$$\rho:\lim_{\substack{\rightarrow \\ M}} H^*(-,\K_{\cdot}^M) \rightarrow H^*(-,\F_{\cdot}).$$
If
\begin{enumerate*}[$\blacktriangleright$]
\item $\rho(i_N)=\psi(\lambda^k(E^N_{\cdot}-N)) $, and, \item
$\rho(u_r)=u_{\lambda^k(r)}$,
\end{enumerate*}
then, $\rho$ agrees with $\lambda_{GS}^k:
K_m(X_{\cdot})\rightarrow K_m(X_{\cdot})$, for every degenerate
simplicial scheme $X_{\cdot}$  in $\mathbf{C}$.
\end{theo}

\begin{theo}[Adams operations]\label{uni3} Let $\{\rho_N:\K_{\cdot}^N\rightarrow \K_{\cdot}\}_{N\geq 1}$ be a weakly additive
system of maps with respect to the $H$-sum of $\K_{\cdot}$. Let
$\rho$ be the induced morphism
$$\lim_{\substack{\rightarrow \\ M}} H^*(-,\K_{\cdot}^M) \rightarrow H^*(-,\F_{\cdot}).$$ If
\begin{enumerate*}[$\blacktriangleright$]
\item $\rho(i_N)=\psi(\Psi^k(E_{\cdot}^N-N)) $, and, \item
$\rho(u_r)=u_{\psi^k(r)}$,
\end{enumerate*}
then, $\rho$ agrees with $\Psi_{GS}^k: K_m(X_{\cdot})\rightarrow
K_m(X_{\cdot})$,  for every   degenerate  simplicial scheme
$X_{\cdot}$   in $\mathbf{C}$.
\end{theo}

Since the Adams operations are group morphisms, it is natural to expect that
they will be induced by $H$-space maps
$$\K_{\cdot}\rightarrow \K_{\cdot} $$
in $\Ho(\mathbf{sT}_*)$. The next two corollaries follow easily from the last
theorem.

\begin{cor}
Let $\rho: \K_{\cdot}\rightarrow \K_{\cdot} $ be an $H$-space map
in the homotopy category of simplicial sheaves on $\mathbf{C}$. If
\begin{enumerate*}[$\blacktriangleright$]
\item $\rho(j_N)=\psi(\Psi^k(E_N-N)) $, and, \item
$\rho(u_r)=\psi(\Psi^k(u_r))$,
\end{enumerate*} then
$\rho$ agrees with the Adams operation $\Psi^k_{GS}$,  for all
  degenerate simplicial schemes  in $\mathbf{C}$.
\end{cor}

Let $S_{\cdot}\mathcal{P}$ denote the Waldhausen simplicial sheaf
on $ZAR(S)$ given by
$$X\mapsto S_{\cdot}\mmP(X)=S_{\cdot}(X).$$

\begin{cor}\label{adams}
Let $\rho: S_{\cdot}\mmP\rightarrow S_{\cdot}\mmP $ be an
$H$-space map in $\Ho(\mathbf{sT}_*)$. If for some $k\geq 1$ there
is a commutative square
$$\xymatrix{  K_0(\Vect(B_{\cdot}GL_{N/S}))\ar[r]^(.55){\psi} \ar[d]_{\Psi^k} &
 K_0(B_{\cdot}GL_{N/S}) \ar[d]^{\rho}  \\
K_0(\Vect(B_{\cdot}GL_{N/S})) \ar[r]_(.55){\psi} &
K_0(B_{\cdot}GL_{N/S}), }$$ then $\rho$ agrees with the Adams
operation $\Psi^k_{GS}$, for all degenerate simplicial schemes
 in $\mathbf{C}$.
\end{cor}

Therefore, there is a unique way to extend the Adams operations from the Grothendieck group of simplicial schemes to higher $K$-theory by means of a
map $S_{\cdot}\mmP\rightarrow S_{\cdot}\mmP$ in the homotopy category of simplicial sheaves.

 Grayson, in \cite{Grayson1}, defines the Adams operations for the
$K$-groups of any exact category with a suitable notion of tensor, symmetric
 and exterior product. The category of vector bundles
  over a scheme satisfies the required conditions, as well as the category
of vector bundles over a simplicial scheme. For every scheme $X$, he constructs
\begin{enumerate*}[$\blacktriangleright$]
\item two $(k-1)$-simplicial sets, $S_{\cdot}\tilde{G}^{(k-1)}(X)$
and $\Sub_k(X)_{\cdot}$, whose diagonals are weakly equivalent to
$S_{\cdot}(X)$, and  \item a $(k-1)$-simplicial map
$\Sub\nolimits_k(X)_{\cdot}\xrightarrow{\Psi^k}
S_{\cdot}\tilde{G}^{(k-1)}(X).$
\end{enumerate*}
His construction is functorial on $X$ and hence induces a map of presheaves.

Grayson has already checked that the operations that he defined induce the
usual ones for the Grothendieck group of a suitable category $\mathcal{P}$.
Therefore, since the conditions of proposition \ref{adams} are fulfilled, we
obtain the following corollary.

\begin{cor} Let $S$ be a  finite dimensional noetherian scheme.
The Adams operations defined by Grayson in \cite{Grayson1} agree
with the Adams operations defined by Gillet and Soul\'e in
\cite{GilletSouleFiltrations}, for every scheme in $\ZAR(S)$. In
particular, they satisfy the usual identities for schemes in
$\ZAR(S)$.
\end{cor}

Grayson did not prove that his operations satisfied the identities
of a lambda ring. It follows from the previous corollary that they
are satisfied for finite dimensional noetherian schemes.

\section{Morphisms between K-theory and cohomology}

\subsection{Sheaf cohomology as a generalized cohomology theory}
Fix $\mathbf{C}$ to be a subsite of the big Zariski site $\ZAR(S)$, as in
section \ref{zariski}.

Consider the Dold-Puppe functor $\mathcal{K}_{\cdot}(\cdot)$
(see \cite{DoldPuppe}), which associates to every
cochain complex of abelian groups concentrated in non-positive
degrees, $G^*$, a simplicial abelian group $\mmK_{\cdot}(G)$,
pointed by zero. It satisfies the property that
$\pi_i(\mmK_{\cdot}(G),0)=H^{-i}(G^*)$.

Now let $G^*$ be an arbitrary cochain complex. Let $(\tau_{\leq
n}G)[n]^*$ be the truncation at degree $n$ of $G^*$ followed by
the translation by $n$. That is,
$$
(\tau_{\leq n}G)[n]^i = \left\{\begin{array}{ll} G^{i+n} & \textrm{if }i<0, \\
\ker (d:G^{n}\rightarrow G^{n+1}) & \textrm{if }i=0, \\
0 & \textrm{if }i>0.
\end{array} \right.
$$
One defines a simplicial abelian group by
$$\mmK_{\cdot}(G)_n:= \mmK_{\cdot}((\tau_{\leq n}G)[n]).$$
The simplicial abelian groups $\mmK_{\cdot}(G)_n$ form an infinite
loop spectrum. Moreover, this construction is functorial on $G$.

Let $\mathcal{F}^*$ be a cochain complex of sheaves of abelian
groups in $\mathbf{C}$, and let $\mmK_{\cdot}(\mmF)_*$ be the
infinite loop spectrum obtained applying section-wise the
construction above. For every $n$, $\mmK_{\cdot}(\mmF)_n$ is an
H-space, since it is a simplicial sheaf of abelian groups.

\begin{lema}[\cite{HuberWild} Prop. B.3.2]
Let $\mmF^*$ be a bounded below complex of sheaves on $\mathbf{C}$
and let $X$ be a scheme in the underlying category. Then, for all
$m\in \Z$,
$$H^{m}(X,\mmK_{\cdot}(\mmF)_*)\cong H^{m}_{\ZAR}(X,\mmF^*).$$
\end{lema}
Here, the right hand side is the usual Zariski cohomology and the
left hand side is the generalized cohomology of the simplicial
sheaf of groups $\mmK_{\cdot}(\mmF)$. Observe that since
$\mmK_{\cdot}(\mmF)$ is an infinite loop space, we can consider
generalized cohomology groups for all integer degrees. Thus we see that the
usual Zariski cohomology can be expressed in terms of generalized
sheaf cohomology using the Dold-Puppe functor.

\subsection{Uniqueness of characteristic classes }\label{uniqueness4}
Now fix a bounded below graded complex of sheaves $\mmF^*(*)$ of
abelian groups, giving a twisted duality cohomology theory in the
sense of Gillet, \cite{Gillet}. In loc. cit., Gillet constructed
\emph{Chern classes} for higher $K$-theory. They are given by a
map of spaces
\begin{equation}\label{cj}
c_j: \K_{\cdot} \rightarrow \mmK_{\cdot}(\mmF(j)[2j]),\qquad j\geq
0.
\end{equation}
 More specifically, they
are given by a map
$$\Z_{\infty}B_{\cdot}GL \rightarrow \mmK_{\cdot}(\mmF(j)[2j])$$
extended trivially over the $\Z$ component of $\K_{\cdot}$. By
example \ref{trivial2}, these maps are weakly additive. In fact,
they are weakly additive also with respect to the $H$-sum of
$\mmK_{\cdot}(\mmF(j)[2j])$ (see remark \ref{trivialhspace}).

For any  space $X_{\cdot}$, the map induced after taking
generalized cohomology,
$$c_j:  K_m(X_{\cdot}) \rightarrow H^{2j-m}_{\ZAR}(X_{\cdot},\mmF^*(j)),\qquad j\geq 0,$$
is called the \emph{$j$-th Chern class}. They are group morphisms
for $m> 0$ but only maps for $m=0$. In this last case, for any
vector bundle $E$ over a scheme $X$, $c_j(E)$ is the standard
$j$-th Chern class taking values in the given cohomology theory.

Using the standard formulas on the Chern classes, one obtains the
\emph{Chern character}
$$\ch :  K_m(X_{\cdot}) \rightarrow \prod_{j\geq 0} H^{2j-m}_{\ZAR}(X_{\cdot},\mmF^*(j))\otimes \Q,$$
which is now a group morphism for all $m\geq 0$. It is induced by
an $H$-space map $$ \ch:\K_{\cdot} \rightarrow \prod_{j\geq
0}\mmK_{\cdot}(\mmF(j)[2j])\otimes \Q.$$ The restriction of $\ch$
to $K_0(X)$ is the usual Chern character of a vector bundle.

We will now state the theorems equivalent to theorems \ref{uni4}
and \ref{uni3}, for maps from $K$-theory to cohomology. In order
to do this, we should first understand better  $c_j(i_N)$ and
$\ch(i_N)$ for all $j,N$. This will be achieved by means of the
Grassmanian schemes.

Denote by
\begin{eqnarray*}
\mmF^*(*) &=& \prod_{i\geq 0,\ j\in \Z}\mmF^i(j), \\
H^*_{\ZAR}(X_{\cdot},\mmF^*(*)) &=& \prod_{i\geq 0,\ j\in \Z}
H^i_{\ZAR}(X_{\cdot},\mmF^*(j)).
\end{eqnarray*}

Let $Gr(N,k)=Gr_{\Z}(N,k)\times_{\Z}S$ be the Grassmanian scheme
over $S$. This is a projective scheme over $S$. Consider $E_{N,k}$
the rank $N$ universal bundle of $Gr(N,k)$ and
$\{U_k\}\xrightarrow{p} Gr(N,k)$ its standard trivialization.
There is a classifying map of the vector bundle $E_{N,k}$,
$\varphi_k: N_{\cdot}U_k\rightarrow B_{\cdot}GL_{N/S}$, satisfying
$p^*(E_{N,k})=\varphi_k^*(E^N_{\cdot})$. This map induces a map in
the Zariski cohomology
$$ H^{*}_{\ZAR}(B_{\cdot}GL_{N/S},\mmF^*(*))\xrightarrow{\varphi_k^*} H^{*}_{\ZAR}(N_{\cdot}U_k,\mmF^*(*))
\cong H^{*}_{\ZAR}(Gr(N,k),\mmF^*(*)).$$ Moreover, for each $m_0$,
there exists $k_0$ such that if $m\leq m_0$ and $k\geq k_0$,
$\varphi_k^*$ is an isomorphism on the $m$-th cohomology group.

\begin{prop}\label{grassmanians}
Let $\chi_1=\{\chi_1^N\}$ and $\chi_2=\{\chi_2^N\}$ be two weakly additive
systems of maps
$$\chi_i^N:\K_{\cdot}^N\rightarrow \mmK_{\cdot}(\mmF(*)),\quad i=1,2, $$
with respect to the same operation. Then, the induced maps
$$\chi_1,\chi_2: K_m(X) \rightarrow H^{*}_{\ZAR}(X,\mmF^*(*)) $$
agree for every scheme $X$ in $\mathbf{C}$, if and only if they agree for
$X=Gr(N,k)$, for all $N$ and $k$.
\end{prop}
\begin{proof} One implication is obvious. For the other implication, fix $m_0$ and let $k_0$ be an integer such that for
every $k\geq k_0$ there is an isomorphism at the $m_0$ level. Then, there are
commutative diagrams
$$\xymatrix{ \lim_{\substack{\rightarrow \\ M}} H^{-m_0}(B_{\cdot}GL_{N/S},\K_{\cdot}^M) \ar[r]^(.52){\chi_1,\chi_2}
\ar[d]_{\varphi^*_k} & H^{*}_{\ZAR}(B_{\cdot}GL_{N/S}, \mmF^*(*))
\ar[d]^{\varphi^*_k}_{\cong}
\\ \lim_{\substack{\rightarrow \\ M}} H^{-m_0}(N_{\cdot}U_k,\K_{\cdot}^M)\ar[r]^(.52){\chi_1,\chi_2}
 & H^{-m_0}_{\ZAR}(N_{\cdot}U_k,\mmF^*(*))  \\ H^{-m_0}(Gr(N,k),\K_{\cdot}) \ar[u]^{p^*}_{\cong}\ar[r]^(.52){\chi_1,\chi_2}
 & H^{-m_0}_{\ZAR}(Gr(N,k),\mmF^*(*)) \ar[u]_{p^*}^{\cong}. }$$
By theorem \ref{uni7},  $\chi_1=\chi_2$ for all schemes $X$, if
they agree for $B_{\cdot}GL_{N/S}$ for all $N\geq 1$. Let $x\in
\lim_{\substack{\rightarrow \\ M}}
H^{-m_0}(B_{\cdot}GL_{N/S},\K_{\cdot}^M)$. Then,
\begin{eqnarray*}
\chi_1(x)=\chi_2(x) & \Leftrightarrow &
(p^*)^{-1}\varphi_k^*\chi_1(x)=(p^*)^{-1}\varphi_k^*\chi_2(x)\\
&\Leftrightarrow &
\chi_1(p^*)^{-1}\varphi_k^*(x)=\chi_2(p^*)^{-1}\varphi_k^*(x),\end{eqnarray*}
and since they agree for the Grassmanians, the proposition is
proved.
\end{proof}

The following two theorems follow from the results \ref{uni7} and
\ref{uni}, together with the preceding proposition.

\begin{theo}[Chern classes]\label{uni5}
There is a unique way to extend the $j$-th Chern class of vector
bundles over schemes in $\mathbf{C}$, by means of a weakly
additive system of maps $\{\rho_N:\K_{\cdot}^N\rightarrow
\mmK_{\cdot}(\mmF(j)[2j])\}_{N\geq 1}$ with respect to the
$H$-space operation in $\mmK_{\cdot}(\mmF(*))$.
\end{theo}

Observe that it follows from the theorem that any weakly additive
collection of maps  with respect to the $H$-space operation of
$\mmK_{\cdot}(\mmF(*))$, inducing the $j$-th Chern class, is
necessarily trivial on the $\Z$-component for $j>0$.

\begin{theo}[Chern character]\label{uni6}
Let $$\K_{\cdot}\longrightarrow \prod_{j\in \Z}
\mmK_{\cdot}(\mmF(j)[2j])$$ be an $H$-space map in
$\Ho(\mathbf{sT}_*)$. The induced morphisms
$$K_m(X)\rightarrow \prod_{j\in \Z} H^{2j-m}_{\ZAR}(X,\mmF^*(j)) $$
agree with the Chern character defined by Gillet in \cite{Gillet} for every
scheme $X$, if and only if, the induced map
$$K_0(X)\rightarrow \prod_{j\in \Z} H^{2j}_{\ZAR}(X,\mmF^*(j)) $$
is the Chern character for $X=Gr(N,k)$, for all $N,k$.
\end{theo}

\begin{cor}
There is a unique way to extend the standard Chern character of vector bundles
over schemes in $\mathbf{C}$, by means of an $H$-space map
$$\rho:\K_{\cdot}\rightarrow \prod_{j\in \Z} \mmK_{\cdot}(\mmF(j)[2j]).$$
\end{cor}

We deduce from these theorems that any simplicial sheaf map
$$S_{\cdot}\mmP \rightarrow \mmK_{\cdot}(\mmF(*)) $$
that induces either the Chern character or any Chern class map at
the level of $K_0(X)$, induces the Chern character or the Chern
class map on the higher $K$-groups of $X$.

\begin{obs} Let  $\mathbf{C}$ be the site of smooth complex varieties and let
 $\mmD^*(*)$ be a graded complex computing absolute Hodge
cohomology. Burgos and Wang, in \cite{Burgos1}, constructed a
simplicial sheaf map $S_{\cdot}\mmP \rightarrow
\mmK_{\cdot}(\mmD(*)) $ which induces the Chern character on any
smooth proper complex variety. A consequence of the last corollary
is that their definition agrees with the \emph{Beilinson
regulator} (the Chern character for absolute Hodge cohomology).

This is not a new result. Using other methods, Burgos and Wang
already proved that the morphism they defined was the same as the
Beilinson regulator. The result is proved there by means of the
bisimplicial scheme $B_{\cdot}P$ introduced by Schechtman in
\cite{Sch}. This introduced an unnecessary delooping, making the
proof generalizable only to sheaf maps inducing group morphisms
and introducing irrelevant ingredients to the proof.
\end{obs}


\begin{thebibliography}{SGA73}

\bibitem[BF08]{FeliuChow}
J.~I. Burgos and E.~Feliu, \emph{Higher arithmetic {C}how groups}, Preprint
  (2008).

\bibitem[BG73]{BrownGestern}
K.~S. Brown and S.~M. Gersten, \emph{Algebraic {$K$}-theory as generalized
  sheaf cohomology}, Algebraic K-theory, I: Higher K-theories (Proc. Conf.,
  Battelle Memorial Inst., Seattle, Wash., 1972), Springer, Berlin, 1973,
  pp.~266--292. Lecture Notes in Math., Vol. 341.

\bibitem[BK72]{bousfieldkan}
A.~K. Bousfield and D.~M. Kan, \emph{Homotopy limits, completions and
  localizations}, Springer-Verlag, Berlin, 1972, Lecture Notes in Mathematics,
  Vol. 304.

\bibitem[BW98]{Burgos1}
J.~I. Burgos and S.~Wang, \emph{Higher {B}ott-{C}hern forms and {B}eilinson's
  regulator}, Invent. Math. \textbf{132} (1998), no.~2, 261--305.

\bibitem[DP61]{DoldPuppe}
A.~Dold and D.~Puppe, \emph{Homologie nicht-additiver {F}unktoren.
  {A}nwendungen}, Ann. Inst. Fourier Grenoble \textbf{11} (1961), 201--312.

\bibitem[Fel07]{FeliuThesis}
E.~Feliu, \emph{On higher arithmetic intersection theory}, PhD Thesis,
  Available at: http://www.tdx.cat/TDX-1220107-112706, 2007.

\bibitem[Fel08]{Feliu-Adams}
\bysame, \emph{{A}dams operations on higher arithmetic k-theory}, Preprint
  (2008).

\bibitem[Gil81]{Gillet}
H.~Gillet, \emph{Riemann-{R}och theorems for higher algebraic {$K$}-theory},
  Adv. in Math. \textbf{40} (1981), no.~3, 203--289.

\bibitem[Gon05]{Goncharov}
A.~B. Goncharov, \emph{Polylogarithms, regulators, and {A}rakelov motivic
  complexes}, J. Amer. Math. Soc. \textbf{18} (2005), no.~1, 1--60
  (electronic).

\bibitem[Gra92]{Grayson1}
D.~R. Grayson, \emph{Adams operations on higher {$K$}-theory}, $K$-Theory
  \textbf{6} (1992), no.~2, 97--111.

\bibitem[GS99]{GilletSouleFiltrations}
H.~Gillet and C.~Soul{\'e}, \emph{Filtrations on higher algebraic
  {$K$}-theory}, Algebraic $K$-theory (Seattle, WA, 1997), Proc. Sympos. Pure
  Math., vol.~67, Amer. Math. Soc., Providence, RI, 1999, pp.~89--148.

\bibitem[Hir03]{Hirschhorn}
P.~S. Hirschhorn, \emph{Model categories and their localizations}, Mathematical
  Surveys and Monographs, vol.~99, American Mathematical Society, Providence,
  RI, 2003.

\bibitem[HW98]{HuberWild}
A.~Huber and J.~Wildeshaus, \emph{Classical motivic polylogarithm according to
  {B}eilinson and {D}eligne}, Doc. Math. \textbf{3} (1998), 27--133
  (electronic).

\bibitem[Jar87]{Jarsimppre}
J.~F. Jardine, \emph{Simplicial presheaves}, J. Pure Appl. Algebra \textbf{47}
  (1987), no.~1, 35--87.

\bibitem[Qui67]{Quillen0}
D.~G. Quillen, \emph{Homotopical algebra}, Lecture Notes in Mathematics, No.
  43, Springer-Verlag, Berlin, 1967.

\bibitem[Qui73]{Quillen}
D.~Quillen, \emph{Higher algebraic {$K$}-theory. {I}}, Algebraic $K$-theory, I:
  Higher $K$-theories (Proc. Conf., Battelle Memorial Inst., Seattle, Wash.,
  1972), Springer, Berlin, 1973, pp.~85--147. Lecture Notes in Math., Vol. 341.

\bibitem[Sch87]{Sch}
V.~V. Schechtman, \emph{On the delooping of {C}hern character and {A}dams
  operations}, $K$-theory, arithmetic and geometry (Moscow, 1984--1986),
  Lecture Notes in Math., vol. 1289, Springer, Berlin, 1987, pp.~265--319.

\bibitem[SGA73]{SGA4}
\emph{Th\'eorie des topos et cohomologie \'etale des sch\'emas. {T}ome 3},
  Springer-Verlag, Berlin, 1973, S\'eminaire de G\'eom\'etrie Alg\'ebrique du
  Bois-Marie 1963--1964 (SGA 4), Dirig\'e par M. Artin, A. Grothendieck et J.
  L. Verdier. Avec la collaboration de P. Deligne et B. Saint-Donat, Lecture
  Notes in Mathematics, Vol. 305.

\bibitem[Sou85]{Soule}
C.~Soul{\'e}, \emph{Op\'erations en {$K$}-th\'eorie alg\'ebrique}, Canad. J.
  Math. \textbf{37} (1985), no.~3, 488--550.

\bibitem[Wal78]{Waldhausen}
F.~Waldhausen, \emph{Algebraic {$K$}-theory of generalized free products. {I},
  {II}}, Ann. of Math. (2) \textbf{108} (1978), no.~1, 135--204.

\end{thebibliography}

\providecommand{\bysame}{\leavevmode\hbox to3em{\hrulefill}\thinspace}
\providecommand{\MR}{\relax\ifhmode\unskip\space\fi MR }
\providecommand{\MRhref}[2]{%
  \href{http://www.ams.org/mathscinet-getitem?mr=#1}{#2}
}
\providecommand{\href}[2]{#2}

\end{document}